\documentclass[reqno]{amsart}
\usepackage{mathrsfs}
\usepackage{amsmath,amsfonts,amssymb,amsthm}
\usepackage{upgreek}
\usepackage[usenames,dvipsnames]{xcolor}
\usepackage{enumitem}
\usepackage{graphicx}
\usepackage{wrapfig}
\usepackage{overpic}
\usepackage{cancel}
\usepackage{float}
\usepackage[outdir=./]{epstopdf}
\usepackage[colorlinks=true]{hyperref} 
\usepackage{cleveref}
\usepackage{url}
\usepackage[utf8]{inputenc}
\usepackage{booktabs,multirow}
\usepackage[font=footnotesize,width=\textwidth]{caption}
\hypersetup{linkcolor=Violet,citecolor=PineGreen}
\newtheorem{teor}{Theorem}[section]

\newtheorem{prop}[teor]{Proposition}

\theoremstyle{definition}

\newtheorem{nota}[teor]{Remark}
\newtheorem{notas}[teor]{Remarks}
\numberwithin{equation}{section}

\newcommand{\R}{\mathbb R}

\newcommand{\Z}{\mathbb{Z}}

\newcommand{\C}{\mathbb{C}}
\newcommand{\N}{\mathbb{N}}

\newcommand{\Y}{\mathbb{Y}}

\newcommand{\mA}{\mathcal{A}}
\newcommand{\mB}{\mathcal{B}}
\newcommand{\mC}{\mathcal{C}}

\newcommand{\mK}{\mathcal{K}}
\newcommand{\mM}{\mathcal{M}}
\newcommand{\mN}{\mathcal{N}}
\newcommand{\mO}{\mathcal{O}}

\newcommand{\mR}{\mathcal{R}}

\newcommand{\mU}{\mathcal{U}}
\newcommand{\mV}{\mathcal{V}}
\newcommand{\ep}{\varepsilon}

\newcommand{\mI}{\mathcal{I}}
\newcommand{\mJ}{\mathcal{J}}

\newcommand{\W}{\Omega}
\newcommand{\WR}{\W\times\R}
\newcommand{\RR}{\R\times\R}
\newcommand{\w}{\omega}

\newcommand{\mb}{\mathfrak{b}}

\newcommand{\ml}{\mathfrak{l}}
\newcommand{\mm}{\mathfrak{m}}
\newcommand{\muk}{\mathfrak{u}}

\newcommand{\upalfa}{$\upalpha$}
\newcommand{\upomeg}{$\upomega$}

\newcommand{\ws}{\w{\cdot}s}
\newcommand{\wt}{\w{\cdot}t}
\newcommand{\bwt}{\bar\w{\cdot}t}

\newcommand{\lsm}{\left[\begin{smallmatrix}}
\newcommand{\rsm}{\end{smallmatrix}\right]}

\begin{document}
\title[Global bifurcation diagrams for third-degree polynomial ODEs]
{Global bifurcation diagrams for coercive third-degree polynomial ordinary differential
equations with recurrent nonautonomous coefficients}
\author[C. Elia]{Cinzia Elia}
\author[R. Fabbri]{Roberta Fabbri}
\author[C. N\'{u}\~{n}ez]{Carmen N\'{u}\~{n}ez}
\address[C. Elia]{Dipartimento di Matematica, Universit\`{a} di Bari, Via Orabona 4, 70125 Bari, Italy}
\address[R. Fabbri]{Dipartimento di Matematica Ulisse Dini, Universit\`{a} degli Studii di Firenze,
        Viale Giovanni Battista Morgagni 67/a, 50134 Firenze, Italy}
\address[C. N\'{u}\~{n}ez]{Departamento de Matem\'{a}tica Aplicada, Universidad de Va\-lladolid,
        Paseo Prado de la Magdalena 3-5, 47011 Valladolid, Spain}
\email[C.~Elia]{cinzia.elia@uniba}
\email[R.~Fabbri]{roberta.fabbri@unifi.it}
\email[C.~N\'{u}\~{n}ez]{carmen.nunez@uva.es}
\thanks{C. Elia partially supported by the GNCS-Indam group and the Ministero dell'Universit\`{a} e delle Ricerca  under research grant PRIN2022PNR P2022M7JZW \ SAFER MESH - "Sustainable mAnagement oF watEr Resources modEls and numerical MetHods". C. N\'{u}\~{n}ez was supported by Ministerio de Ciencia, Innovaci\'{o}n y Universidades
(Spain) under project PID2021-125446NB-I00 and by Universidad de Valladolid under project PIP-TCESC-2020, and
she is a member of IMUVA.}
\date{}
\begin{abstract}
Nonautonomous bifurcation theory is a growing branch of mathematics, for the insight
it provides into radical changes in the global dynamics of realistic models for many real-world
phenomena, i.e., into the occurrence of critical transitions.
This paper describes several global bifurcation diagrams for nonautonomous first order
scalar ordinary differential equations generated by coercive third degree
polynomials in the state variable. The conclusions are applied to a population dynamics model
subject to an Allee effect that is weak in the absence of migration and becomes strong under a
migratory phenomenon whose sense and intensity depend on a threshold in the number of
individuals in the population.
\end{abstract}
\renewcommand{\subjclassname}{\textup{2020} Mathematics Subject Classification}
\maketitle
\section{Introduction}\label{sec:introduccion}
Bifurcation theory is a branch of the study of dynamical systems that dates back
to the early works of Poincar\'{e} \cite{poin} at the end of the 19th century. Much more recent is
the extension of this theory to non-autonomous dynamical systems. The analysis of
these systems arises from the need of the applied branches of science to describe
models whose own laws of evolution change with respect to time, which generally
allows a more realistic description of the phenomenon.  All these models depend
on parameters, and very often a small variation in one of these parameters causes
a strong variation in the resulting global dynamics. Understanding the mechanisms
of occurrence of these changes and, closely related, describing the dynamics for
close values of the parameter are, broadly speaking, the objectives of
nonautonomous bifurcation theory.

The most common approach to autonomous bifurcation theory for one-parametric families of
{\em scalar\/} ordinary differential equations (ODEs in what follows) analyzes the evolution,
as the parameter varies, of the number and type of critical points, which correspond to the
constant solutions of the ODEs. They are classified into hyperbolic attractive, hyperbolic repulsive,
and nonhyperbolic, and often determine the global phase line. The object of study is not so clear in the
nonautonomous extension of the theory, since a scalar time-dependent ODE $x'=f(t,x)$ does not admit,
in general, constant solutions. So, although the overall objective is basically always the same, there is not
total agreement on where to place the focus for the analysis. Different approaches
are presented in the works of Braaksma et al.~\cite{brbh},
Johnson and Mantellini \cite{joma}, Fabbri et al.~\cite{fajm}, Kloeden \cite{kloe},
Langa et al.~\cite{lars}, Rasmussen \cite{rasm1,rasm2}, N\'{u}\~{n}ez and Obaya \cite{nuob8,nuob10},
J\"{a}ger \cite{jage3}, P\"{o}tzsche \cite{potz2,potz4}, Kloeden and Rasmussen \cite{klra},
Anagnostopoulou and J\"{a}ger \cite{anja},
Anagnostopoulu et al. \cite{anjk, anpr}, Fuhrmann \cite{fuhr}, Longo et al.~\cite{lnor},
Remo et al. \cite{refj}, and Due\~{n}as et al.~\cite{duno1,duno2,duno5}, as well as in the references therein.

In this work, following in the wake of \cite{duno1,duno2,duno5}, we analyze the bifurcation problem given
by the variation in $\ep$ of an $\ep$-parametric family of third degree coercive polynomial nonautonomous ODEs,
\begin{equation}\label{eq:intro}
 x'=-x^3+\bar c(t)\,x^2+\ep\,\big(\bar b(t)\,x+\bar a(t)\big)\,,
\end{equation}
determined by three bounded and uniformly continuous maps $\bar c,\,\bar b,\,\bar a\colon\R\to\R$.
With the approach previously established in \cite{nuob8,nuob10,lnor},
we use the skew-product formalism, defining from \eqref{eq:intro} a (possibly local) real continuous
flow $\tau_\ep$ on the trivial bundle $\WR$, where $\W$ is the hull of $(\bar c,\,\bar b,\,\bar a)$.
That is, $\W$ is the (compact) closure in the compact-open topology of $C(\R,\R^3)$ of the set of time-shifts
$\{(\bar c,\bar b,\bar a){\cdot}t\,|\;t\in\R\}$, where $\bar d{\cdot}t (s)=\bar d(t+s)$.
Defining $c(\w)=\w_1(0)$, $b(\w)=\w_2(0)$ and $a(\w)=\w_3(0)$ for $\w=(\w_1,\w_2,\w_3)\in\W$,
we obtain, for each $\ep\in\R$, the family of equations
\begin{equation}\label{eq:introhull}
 x'=-x^3+c(\wt)\,x^2+\ep\,\big(b(\wt)\,x+a(\wt)\big)\,,\qquad\w\in\W\,,
\end{equation}
whose solutions $v_\ep(t,\w,x)$ satisfying $v_\ep(0,\w,x)=x$ yield the fiber-component of the
flow $\tau_\ep$, which is of skew-product type: $\tau_\ep(t,\w,x)=(\wt,v_\ep(t,\w,x))$.
Observe that \eqref{eq:intro} is \eqref{eq:introhull} for $\w=(\bar c,\bar b,\bar a)$.
We assume the time-shift flow on the hull $\W$ to be minimal and uniquely ergodic
(which is the situation in many nonautonomous mathematical models, as those determined by an
almost periodic function $(\bar c,\bar b,\bar a)\colon\R\to\R^3$),
and choose the minimal subsets of $\WR$ as the objects whose variation in number and type
(hyperbolic attractive,
hyperbolic repulsive of nonhyperbolic) determine the occurrence of bifurcation values of $\ep$.
The most basic minimal set is the so-called {\em copy of the base}, which is the (invariant) graph of
a continuous map $\mb_\ep\colon\W\to\R$ such that $\mb_\ep(\wt)=v_\ep(t,\w,\mb_\ep(\w))$ for all $(t,\w)\in\R\times\W$: this is the natural extension of a critical value in the autonomous case.
So, our approach is quite natural, although unlike in the autonomous framework there may be
minimal subsets with a highly complicated dynamics. Some of the first samples of these extremely complex
minimal sets, that include strange non-chaotic attractors,
can be found in Million\u{s}\u{c}ikov~\cite{mill,mill2}, Vinograd \cite{vino}
(see also Lipnitskii \cite{lipn} for some technical improvements),
Johnson~\cite{john12}, and Koltyzhenkov~\cite{kolt} (and in Grebogi {\em et al.} \cite{gopy},
Bezhaeva and Oseledets \cite{beos}, and Keller~\cite{kell} for discrete instead of continuous flows).
In the bifurcation diagrams described in
\cite{nuob8,nuob10,lnor,duno1,duno2,duno5} we observe a phenomenon that appears frequently in the
literature: these complex sets can appear only at the bifurcation values of the parameter. This will
also be the situation of the problem studied here. So, we get one more sample that the degree of
rarity of these sets depends not only on their intern dynamics, but also on their extreme lack
of persistence under quite standard perturbations.
\par
Returning to the particular case of \eqref{eq:introhull},
it is enough to work with autonomous examples
to see that the possibilities of the bifurcation diagram are very numerous.
In order for the result of our analysis to be of reasonable length, we need to make certain
choices at the beginning. Modifying these assumptions will substantially change the results,
but the study of many of the cases that we do not consider in this paper can be done using the same techniques:
classical general methods of topological dynamics and ergodic theory combined with new results and
techniques, in the line of those developed in \cite{nuob8, nuos4, duno1, duno5}.
The main results of this paper are obtained under the conditions $\inf_{\w\in\W}d(\w)>0$ for
$d=c,\,b,\,-a$; and, like in the autonomous case (with $c>0$ $b>0$ and $a<0$),
the relative sizes of $c$ and $-a/b$ determine very different bifurcation situations.

In all these autonomous bifurcation diagrams, only two types of bifurcations appear:
local saddle-node bifurcations, when two branches of hyperbolic critical points exist to
the left (or right) of $\ep_0$ and collide at this point, giving rise to a unique
nonhyperbolic critical point at $\ep_0$ and to the local absence of critical points to its right (or left);
and local transcritical bifurcation, when two branches of hyperbolic critical points exist
both at the left and right of $\ep_0$, and they collide at a unique nonhyperbolic critical point at $\ep_0$.
In the nonautonomous setting, we say that $\ep_0$ is a local saddle-node bifurcation point
when two hyperbolic copies of the base which exist for close values of $\ep<\ep_0$ (or $\ep>\ep_0)$ approach  each
other as $\ep\to(\ep_0)^-$ (or as $\ep\to(\ep_0)^+$) until they collide at least at a point,
giving rise to a locally unique $\tau_{\ep_0}$-minimal set, which is
nonhyperbolic, and to the absence of minimal sets ``nearby" for close $\ep>\ep_0$ (or $\ep<\ep_0)$.
And we say that $\ep_0$ is a local transcritical bifurcation point
when two hyperbolic copies of the base exist for close values of $\ep$  and approach each
other as $\ep\to\ep_0$ until they collide at $\ep_0$, giving rise to a locally unique
$\tau_{\ep_0}$-minimal set which is nonhyperbolic. For our problem, roughly speaking, we prove that
\begin{itemize}[leftmargin=10pt]
\item[-] $0$ is always a local saddle-node bifurcation point, which appears as the result
of the global collision of two hyperbolic copies of the base as $\ep\to(0)^+$.
\item[-] When $\sup_{\w\in\W}c(\w)<\inf_{\w\in\W}(-a(\w)/b(\w))$, there are at least two more values
of (possibly partial) collision of hyperbolic copies of the base: $\ep^*>\ep_*>0$; the three values are
local saddle-node bifurcation points if, in addition, $a$ is a real multiple of $b$; and they are the
unique ones if the oscillation of $c$ is not too strong.
\item[-] When $\inf_{\w\in\W}c(\w)>\sup_{\w\in\W}(-a(\w)/b(\w))$, there are no strictly positive bifurcation values.
Additional conditions determine either the absence of
negative bifurcation values or the existence of exactly two of them, also of saddle-node type.
\item[-] When $c(\w)=-a(\w)/b(\w)=s$, with $s$ constant, there is a strictly positive bifurcation value,
of local transcritical type, and none negative.
\end{itemize}
The results outlined above are better understood by having a look to the depictions in Figures \ref{fig:caso1}, \ref{fig:casos23} and \ref{fig:caso4}, in Section \ref{sec:results}.
In all the situations, the analysis also involves a description of the evolution of the global
attractor $\mA_\ep$ of the flow $\tau_\ep$ and, in most of the cases, the bifurcation values are points
of discontinuity of the map $\ep\mapsto\mA_\ep$. The hypotheses and results, assumed and proved
for the skew-product, are easily rewritten for the initial family \eqref{eq:intro}. In this
reformulation, instead of considering the evolution in the type and
number of minimal sets, we focus on the number and type of hyperbolic solutions.

The results are applied to describe the evolution of a single population
in a given habitat, subject to an Allee effect (see Courchamp at al.~\cite{cobg})
which is weak in the absence of migration,
and to a particular type of migration whose intensity depends on a threshold
in the number of individuals in the habitat.
The bifurcation points can be read in terms of critical transitions (see~Scheffer \cite{scheffer3}):
significant changes in the state of a complex system that occur as consequences of
small variations in its inputs.

We complete the Introduction with a brief sketch of the structure of the paper.
Section \ref{sec:preliminares} contains the basic concepts and properties required to understand the rest of the
paper: we introduce the skew-product framework we work in; we recall the concepts of equilibria, hyperbolic and nonhyperbolic
minimal set, and global attractor; we summarize some properties of the Lyapunov exponents; and we
describe with more detail the hull construction outlined above. The core of the paper is
Section \ref{sec:results}, where we obtain the global bifurcation diagrams mentioned above
(and additional results under less restrictive hypotheses, not described in these first paragraphs),
and where we indicate how to particularize each of them to a parametric family of processes instead of flows.
Finally, in Section \ref{sec:model}, we apply our previous results to analyze the occurrence of
critical transitions in a particular population dynamics model.
\section{Some preliminary results}\label{sec:preliminares}
In this section we recall the main concepts and tools required to prove the main results
in Section \ref{sec:results}. A {\em (real and continuous) flow\/} on a topological space $\Y$ is
a continuous map $\sigma\colon\mV\subseteq\R\times\Y\to\Y$ defined in an open subset $\mV\supseteq\{0\}\times\Y$
such that, for all $y\in\Y$, $\sigma(0,y)=y$ and $\sigma(t+s,y)=\sigma(t,\sigma(s,y))$ if the right-hand
term is defined. It is {\em global\/} if $\mV=\R\times\Y$.
The definitions of orbit, forward and
backward semiorbit, invariant set, \upalfa-limit set and \upomeg-limit set,
which we omit, can be found in the basic texts of topological dynamics, as \cite{elli}.
We also omit the definitions of (regular) invariant and ergodic measures for the flow:
see, e.g., \cite{mane}.
\par
Let $(\W,\sigma)$ be a global flow on a compact metric space $\W$,
and let us consider the family of equations
\begin{equation}\label{eq:2equ}
 x'=f(\wt,x)\,,\qquad\w\in\W\,,
\end{equation}
where $f\in C^{0,1}(\WR,\R)$; i.e., we assume that the partial derivative $f_x$ globally exists,
and that $f$ and $f_x$ are jointly continuous. (In Section \ref{subsec:hull} we will briefly explain
how such a family arises from a suitable single ODE.)
We represent by $\tau$ the (possibly local) skew-product flow induced by \eqref{eq:2equ} on $\WR$, namely
\begin{equation}\label{def:flujo}
 \tau\colon\mV\subseteq\R\times\WR\to\WR\,,\quad(t,\w,x)\mapsto (\wt,v(t,\w,x))
\end{equation}
where $\mV\supseteq\{0\}\times\W\times\R$ is an open subset and $v(t,\w,x)$
is the maximal solution of the equation \eqref{eq:2equ}
corresponding to $\w$ with initial condition $v(0,\w,x)=x$. We will write $v'(t,\w,x)=f(\wt,v(t,\w,x))$. So,
$v'$ represents $(d/dt)\,v$.
\subsection{Compact invariant sets, upper and lower solutions, and global attractor}\label{subsec:minimal}
The next concepts will play a fundamental role in
some of the proofs.
A map $\mb\colon\W\to\R$ is a $\tau$-{\em equilibrium} if $\mb(\wt)=v(t,\w,\mb(\w))$ for all
$\w\in\W$ and $t\in\R$.
A {\em $\tau$-copy of the base} is the graph of a continuous $\tau$-equilibrium.
(The prefix $\tau$ will be sometimes omitted if there is no risk of confusion.)

Let $\mK\subset\WR$ be a compact $\tau$-invariant set projecting onto the whole set $\W$.
The set $\mK$ is {\em pinched} if the section
$(\mK)_\w:=\{x\in\R\,|\;(\w,x)\in\mK\}$ reduces to a point at least at a point $\w\in\W$.
It is easy to check that its {\em lower\/} (resp.~{\em upper}) {\em equilibrium}, given by
$\ml_\mK(\w):=\sup{x\in\R\,|\;(\w,x)\in\mK})$
(resp.~$\muk_\mK(\w):=\inf{x\in\R\,|\;(\w,x)\in\mK}$) is a lower (resp.~upper)
semicontinuous equilibrium, and hence it is continuous at the points of a residual subset of $\W$.

The flow $(\W,\sigma)$ (or the set $\W$) is {\em minimal\/} if every $\sigma$-orbit is dense in $\W$.
A $\tau$-invariant compact subset $\mM\subset\WR$ is $\tau$-minimal if $(\mM,\tau|_{\mM})$ is minimal;
or, equivalently, if the $\tau$-orbit of any element of $\mM$ is dense in $\mM$. Any
$\tau$-invariant compact set contains a $\tau$-minimal set. If $\W$ is minimal, then any
$\tau$-invariant compact set projects on the whole set $\W$, and any copy of the base is a
$\tau$-minimal set. As already mentioned, the copies of the base are the simplest minimal sets,
playing in many cases
the equivalent role of the equilibrium points for autonomous ODEs. We represent by
$\mM=\{\mb\}$ the $\tau$-minimal set defined by a continuous copy of the base $\mb\colon\W\to\R$.

The next result, for a minimal base flow, is basically proved
in \cite[Section 2]{cano}. A more detailed proof is given in
\cite[Proposition 1.32 and Corollary 1.33]{duen2}.
\begin{prop}\label{prop:defM}
Let the flow $(\W,\sigma)$ be minimal.
\begin{itemize}[leftmargin=20pt]
\item[\rm(i)] Let $\mb\colon\Omega\rightarrow\mathbb{R}$ be a
semicontinuous equilibrium and let $\w_0$ be any continuity point of $\mb$. Then,
\begin{equation}\label{eq:precollision}
 \mM^b=\mathrm{cl}_{\WR}\{(\w_0{\cdot}t,\mb(\w_0{\cdot}t))\,|\;t\in\R\}
\end{equation}
is a minimal set, it is independent of the choice of $\w_0$, and its section $(\mM^b)_\w$ reduces to
$\mb(\w)$ for all the points of the residual $\sigma$-invariant subset of $\W$
given by the continuity points of $\mb$. In addition, the sections $(\mM)_\w$ of any $\tau$-minimal set
$\mM\subset\WR$ are singletons for all the points $\w$ in a residual $\sigma$-invariant subset of $\W$.
\item[\rm(ii)] Two different $\tau$-minimal sets $\mM_1$ and $\mM_2$ are fiber-ordered: if $x^0_1<x^0_2$ for
two points $(\w^0,x^0_1)\in\mM_1$ and $(\w^0,x^0_2)\in\mM_2$, then $x_1<x_2$ whenever
$(\w,x_1)\in\mM_1$ and $(\w,x_2)\in\mM_2$.
\end{itemize}
\end{prop}

A {\em bounded global lower solution\/} for $x'=f(\wt,x)$ is a bounded map
$\mb\colon\W\to\R$ such that $t\mapsto\mb(\wt)$ is $C^1$ and $\mb'(\w)\le f(\w,\mb(\w))$ for all
$\w\in\W$, where $\mb'(\w)=(d/dt)\,\mb(\wt)|_{t=0}$, and it is {\em strict} if the inequality is
strict for all $\w\in\W$. By changing the sign of the inequalities we obtain the definition of
({\em strict}) {\em bounded global upper solution\/}.
%

The constant lower and upper solutions $r$ (respectively characterized by the conditions
$f(\w,r)>0$ and $f(\w,r)<0$ for all $\w\in\W$) will be a useful tool for many points in the proofs
of the main results. The next property will be used often:
\begin{prop}\label{prop:encaje}
Let $m_1<m_2$ be real constants, and assume that one of them is a global upper solution and the other one a lower
global solution. Then, there exists a minimal set contained in $\W\times[m_1,m_2]$. If, in addition,
$m_1$ (resp.~$m_2$) is strict, then the minimal set is contained in $\W\times(m_1,m_2]$ (resp.~$\W\times[m_1,m_2)$).
\end{prop}
\begin{proof}
We choose any $\w\in\W$. It is easy to check that $\W\times[m_1,m_2]$ contains the forward
(resp.~backward) semiorbit of $(\w,(m_1+m_2)/2)$ if $m_2$ (resp.~$m_1$) is the global upper solution.
Hence, $\W\times[m_1,m_2]$ also contains a minimal subset of the \upomeg-limit set (resp.~\upalfa-limit set) of
this orbit.
Let us prove the last assertion in the case that $m_2$ is a strict global upper solution,
assuming for contradiction the existence of a point $(\bar\w,m_2)$ in the
\upomeg-limit set of $(\w,(m_1+m_2)/2)$. Then, since $v'(0,\bar\w,m_2)=f(\bar\w,v(0,\bar\w,m_2))=
f(\bar \w,m_2)<(m_2)'=0$,
we get $v(t,\bar\w,m_2)>m_2$ for small $t<0$. But this is impossible,
since the \upomeg-limit set is $\tau$-invariant and
contained in $\W\times[m_1,m_2]$. The proofs of the three remaining cases are analogous.
\end{proof}
A set $\mA\subset\WR$ is the {\em global attractor for the flow $\tau$} if it is a compact $\tau$-invariant set
that attracts every bounded set $\mC\subset\WR$. This attraction property means that all the forward $\tau$-semiorbits of points of $\mC$ are globally defined (i.e., $[0,\infty)\times\mC\subset\mV$) and that $\lim_{t\rightarrow\infty} \text{dist}(\tau_t(\mC),\mA)=0$, where $\tau_t(\mC):=\{\tau(t,\w,x)\,|\;(\w,x)\in\mC\}$ and
\[
 \text{dist}(\mC_1,\mC_2)=\sup_{(\w_1,x_1)\in\mC_1}
 \Big(\inf_{(\w_2,x_2)\in\C_2}\big(\text{\rm dist}_{\WR}((\w_1,x_1),(\w_2,x_2))\big)\Big)\,.
\]
The next properties are proved in \cite[Theorem 5.1]{duno1}.
\begin{teor}\label{th:atractor}
Assume the coercivity condition $\lim_{x\to\pm\infty}f(\w,x)=\pm\infty$ uniformly on $\W$. Then
all the forward semiorbits are global (i.e., $[0,\infty)\times\W\times\R\subset\mV$),
and there exists the global attractor $\mA$ for $\tau$, which is given by the union of the graphs of all the
bounded solutions of the family of equations \eqref{eq:2equ} and takes the form
\[
 \mA=\bigcup_{\w\in\W}\big(\{\w\}\times[\ml(\w), \muk(\w)]\big)\,.
\]
\end{teor}
\begin{notas}\label{rm:comparison}
In what follows we give some properties of global attractors, bounded global upper and lower
solutions and copies of the base that will be needed in Section \ref{sec:results}.
The coercivity condition of Theorem \ref{th:atractor} is assumed.
\par
1.~If there exists $r_1\in\R$ (resp.~$r_2\in\R$)
such that $f(\w,x)>0$ if $x\le r_1$ (resp.~$f(\w,x)<0$ if $x\ge r_2$) for all $\w\in\W$, then
$\mA\subset\W\times(r_1,\infty)$ (resp.~$\mA\subset\W\times(-\infty,r_2)$). Let us prove
the first assertion, assuming for contradiction that $l:=\inf_{\w\in\W}\ml(\w)=\ml(\bar\w)\le r_1$.
Then, $\ml(\bwt)=v(t,\bar\w,\ml(\bar\w))=v(t,\bar\w,l)<l$ for $t<0$, and this is impossible.
Similarly, if $f(\w,x)>0$ for $x< r_1$ and for all $\w\in\W$, then
$\mA\subset\W\times[r_1,\infty)$. The remaining proofs are analogous.
\par
2.~If $\mb\colon\W\to\R$ is $C^1$ along the base orbits and $\mb'(\w)\le f(\w,\mb(\w))$
(resp.~$\mb'(\w)\ge f(\w,\mb(\w))$) for all $\w\in\W$,
then $\mb\le\muk$ (resp.~$\mb\ge\ml$). If $\mb'(\w)<f(\w,\mb(\w))$
(resp.~$\mb'(\w)>f(\w,\mb(\w))$) for all $\w\in\W$, then $\mb<\muk$ (resp.~$\mb>\ml$).
These properties, based on classical comparison arguments, are proved in
\cite[Theorem 5.1(iii)]{duno1}.
\par
3.~If $\mb\colon\W\to\R$ is upper semicontinuous and a bounded global (strict) upper solution,
then its graph is (strictly) above the \upomeg-limit set $\mO$ of any point $(\w_0,\mb(\w_0))$ (i.e., of
the corresponding orbit);
that is, $x\le\mb(\w)$ ($x<\mb(\w)$) for any point $(\w,x)\in\mO$.
And if, in addition, there exists the \upalfa-limit set of a point
$(\w_0,\mb(\w_0))$, then this set is (strictly) above the graph of $\mb$.
Consequently, in the strict case, no point $(\w_0,\mb(\w_0))$ belongs to any
minimal set. Analogous properties
with reverse orders hold in the case of a bounded global (strict) lower solution given by
a lower semicontinuous maps.
The proofs of these properties are based on comparison results: we first prove the
non-strict inequalities, and then deduce the strict ones in the strict cases by
easy contradiction arguments, solving the equation in the reverse sense of the time.
\end{notas}
\subsection{Lyapunov exponents and hyperbolic minimal sets}\label{subsec:Lyap}
Let $\mK\subset\WR$ be a $\tau$-invariant compact set projecting onto the whole set $\W$.
A value $\gamma\in\R$ is a {\em Lyapunov exponent of $\mK$} if there exists
$(\w,x)\in\mK$ such that
\begin{equation}\label{def:2gamma}
 \gamma=\lim_{t\to\pm\infty}(1/t)\int_0^t f_x(\tau(r,\w,x))\,dr\,.
\end{equation}
Let us assume that $(\W,\sigma)$ is uniquely ergodic, and let us call $m$ the unique $\sigma$-invariant
(and ergodic) measure. Using Riesz' Representation Theorem, Kryloff-Bogoliuboff’s Theorem,
Birkhoff's Ergodic Theorem, \cite[Theorem 4.1]{furs} and \cite[Theorem~1.8.4]{arnol},
it is possible to check that $\gamma$ is a Lyapunov exponent of $\mK$ if and only if there
exists an $m$-measurable equilibrium $\mb\colon\W\to\R$ with graph contained in $\mK$ such that
\[
 \gamma=\int_{\W}f_x(\w,\mb(\w))\,dm\,.
\]
A detailed proof of this assertion, in a much more
general case, can be found in \cite[Sections 1.1.3 and 1.2.4]{duen2}.

A $\tau$-copy of the base $\{\mb\}$ is {\em hyperbolic attractive\/} if it is uniformly
exponentially stable (on the fiber) as time increases; i.e., if there exists $\rho>0$, $k\ge 1$
and $\gamma>0$ such that: if, for any $\w\in\W$, $|\mb(\w)-x|<\rho$, then $v(t,\w,x)$ is
defined for all $t\ge 0$, and in addition $|\mb(\wt)-v(t,\w,x)|\le k\,e^{-\gamma\,t}\,|\mb(\w)-x|$ for $t\ge 0$.
Changing $t\ge 0$ by $t\le 0$ and $\gamma$ by $-\gamma$
provides the definition of {\em repulsive hyperbolic} $\tau$-copy of the base.
We will also say that a $\tau$-minimal set is {\em hyperbolic attractive\/} (resp.~{\em repulsive})
if it is a hyperbolic attractive (resp.~repulsive) $\tau$-copy of the base; and, otherwise,
it is {\em nonhyperbolic}.
\begin{nota}\label{rm:comp2}
An attractive (resp.~repulsive) hyperbolic copy of the base $\{\mb\}$ does not intersect the \upalfa-limit set
(resp.~\upomeg-limit set) of any $(\w,x)$ with $x\ne\mb(\w)$. This intuitive property
is proved in \cite[Proposition 2.6(ii)]{duno4}.
\end{nota}
The next result, which will be repeatedly used, is basically proved in
\cite[Corollary 2.10 and Theorem 3.4]{cano}, and a more detailed proof of (i) and (ii)
is included in the proof of \cite[Theorem 1.40]{duen2}.
\begin{teor}\label{th:caract-hiper}
Let the flow $(\W,\sigma)$ be minimal. Then,
\begin{itemize}[leftmargin=20pt]
\item[\rm(i)] a minimal set is hyperbolic attractive if and only if its upper Lyapunov
exponent is negative.
\item[\rm(ii)] A minimal set is hyperbolic repulsive if and only if its lower Lyapunov
exponent is positive.
\item[\rm(iii)] If the coercivity condition of Theorem {\rm\ref{th:atractor}}
holds, then the global attractor $\mA$ is an attractive hyperbolic $\tau$-copy of the base if and
only if all the $\tau$-minimal sets are hyperbolic attractive.
\end{itemize}
In particular, in the uniquely ergodic case, with ergodic measure $m$, a $\tau$-copy of the base
$\{\mb\}$ is: hyperbolic attractive if and only if $\int_\W f_x(\w,\mb(\w))\,dm<0$, hyperbolic
repulsive if and only if $\int_\W f_x(\w,\mb(\w))\,dm>0$, and (hence) nonhyperbolic
if and only if $\int_\W f_x(\w,\mb(\w))\,dm=0$. And, in the conditions of {\rm(iii)},
$\mA$ is an attractive hyperbolic copy of the base if and only if
$\int_\W f_x(\w,\mb(\w))\,dm<0$ for any $m$-measurable bounded $\tau$-equilibrium.
\end{teor}
\subsection{The hull construction}\label{subsec:hull}
Let us now consider a single ODE,
\begin{equation}\label{eq:ODE}
 x'=\bar f(t,x)\,
\end{equation}
where $\bar f\colon\R\times\R\to\R$ belongs to $C^{0,1}(\R,\R)$; i.e., the derivative $\bar f_x$ with respect
to $x$ globally exists, and the restrictions of the maps $f$ and $f_x$ to $\R\times\mJ$ are bounded
and uniformly continuous for any compact set $\mJ\subset\R$.
Let us define $\bar f_t(s,x):=\bar f(t+s,x)$.
The {\em hull\/} $\W$ of $\bar f$ is the closure of the set
$\{\bar f_t\,|\;t\in\R\}$ on the set $C(\RR,\R)$ provided with the compact-open
topology. Then: the set $\W$ is a compact metric space contained in $C^{0,1}(\RR,\R)$,
the time-shift map $\sigma\colon\R\times\W\to\W,(t,\w)\mapsto\w{\cdot}t:=\wt$
defines a global continuous flow, and the map $f$ given by $f(\w,x)=\w(0,x)$
belongs to $C^{0,1}(\WR,\R)$. The proof of these properties can be found in
\cite[Theorem IV.3]{sell} and \cite[Theorem~I.3.1]{shyi}. Note that
\eqref{eq:ODE} is one of the equations of
the corresponding family \eqref{eq:2equ}: it is given by the element $\w=\bar f\in\W$.
Note also that $(\W,\sigma)$ is a {\em transitive flow},
i.e., there exists a dense $\sigma$-orbit: that of the point $\bar f$.
The map $\bar f$ is {\em recurrent} if $(\W,\sigma)$ is a {\em minimal flow}.
\par
The flow $\tau$ given by \eqref{def:flujo} from the family \eqref{eq:2equ} constructed
from \eqref{eq:ODE} is the {\em skew-product flow
induced by $\bar f$ on its hull}. A standard procedure in nonautonomous dynamics is:
to construct this skew-product flow, use techniques from topological dynamics and ergodic theory
to describe the behavior of its orbits, and derive consequences for the dynamics induced by \eqref{eq:ODE}.
This is basically the approach of this paper: the results are formulated for minimal and
uniquely ergodic flows; but then we show how to extract conclusions for a single recurrent
equation giving rise to a uniquely ergodic hull, and apply them to the analysis of a particular model.
\section{Some global bifurcation diagrams}\label{sec:results}
Let $(\W,\sigma)$ be a global real continuous flow on a compact metric space,
minimal and uniquely ergodic, and
let $m$ be the unique $\sigma$-invariant measure on $\W$. Let
$a,b,c\colon\W\to\R$ be continuous maps, and let us consider the
one-parameter family of families of scalar ODEs
\begin{equation}\label{eq:EDOpol}
 x'=p_\ep(\wt,x)\,,\qquad\w\in\W\,,
\end{equation}
where $\ep$ varies in $\R$ and
\begin{equation}\label{def:p}
 p_\ep(\w,x):=-x^3+c(\w)\,x^2+\ep\,\big(b(\w)\,x+a(\w)\big)\,.
\end{equation}
Recall that a family of this type appears by the
hull procedure from a single ODE: see Sections \ref{subsec:hull} and
\ref{subsec:process}. We will write \eqref{eq:EDOpol}$_\ep$ and \eqref{eq:EDOpol}$_\ep^\w$
to refer to the $\w$-family for a fixed $\ep$ and to a particular equation,
respectively. We also represent by $\tau_\ep$ the
(possibly local) skew-product flow induced by \eqref{eq:EDOpol}$_\ep$ on $\WR$, so that
\[
 \tau_\ep\colon\mV_\ep\subseteq\R\times\WR\to\WR\,,\quad(t,\w,x)\mapsto (\wt,v_\ep(t,\w,x))
\]
where $\mV_\ep\supset\{0\}\times\W\times\R$ is open. Note that $p_\ep$ satisfies the most restrictive conditions of \cite{duno1}, which are the coercivity property $\lim_{x\to\pm\infty} p_\ep(\w,x)/x=-\infty$ uniformly on $\W$, and the strict concavity of the derivative of $x\mapsto p_\ep(\w,x)$ for all $\w\in\W$. Some of the results
of that paper, as well as some of \cite{duno4} (in turn, strongly based on \cite{nuos4} and \cite{lnor}),
will be used in the description of the possibilities for the global $\tau_\ep$-dynamics.
As stated in Theorem \ref{th:atractor}, the coercivity condition ensures that $v_\ep(t,\w,x)$ is defined and
bounded for all $t\ge 0$ (i.e., $[0,\infty)\times\W\times\R\subset\mV_\ep$), and there exists
the global attractor $\mA_\ep$ for $\tau_\ep$, which is given by the union of the graphs of all the
bounded solutions of the family of equations \eqref{eq:EDOpol}$_\ep$, and takes the form
\[
 \mA_\ep=\bigcup_{\w\in\W}\big(\{\w\}\times[\ml_\ep(\w), \muk_\ep(\w)]\big)\,.
\]
Recall also (see Section \ref{subsec:minimal}) that the maps $\ml_\ep,\,\muk_\ep\colon\W\to\R$
are lower and upper semicontinuous equilibria, respectively, and this each of them
is continuous at the points of a residual subset of $\W$.
\begin{teor}\label{th:nummin}
There are three possibilities for the number of $\tau_\ep$-minimal sets:
\begin{itemize}[leftmargin=20pt]
\item[\rm(1)] There are exactly three $\tau_\ep$-minimal sets. In this case,
    they are copies of the base: $\{\ml_\ep\}$, $\{\mm_\ep\}$ and $\{\muk_\ep\}$,
    with $\ml_\ep<\mm_\ep<\muk_\ep$. In addition, $\{\ml_\ep\}$ and $\{\muk_\ep\}$
    are hyperbolic attractive and $\{\mm_\ep\}$ is hyperbolic repulsive.
\item[\rm(2)] There are exactly two $\tau_\ep$-minimal sets. In this case,
    there are two possibilities: either $\{\ml_\ep\}$ is an attractive hyperbolic $\tau_\ep$-copy
    of the base and the other one, nonhyperbolic, is constructed as the closure of
    $\{(\w_0{\cdot}t,\muk_\ep(\w_0{\cdot}t))\,|\;t\in\R\}$ for a continuity point $\w_0$ of $\muk_\ep$,
    and it is a pinched set; or $\{\muk_\ep\}$ is an attractive hyperbolic $\tau_\ep$-copy of the
    base and the other one, nonhyperbolic, is constructed as the closure of
    $\{(\w_0{\cdot}t,\ml_\ep(\w_0{\cdot}t))\,|\;t\in\R\}$ for a continuity point $\w_0$ of $\ml_\ep$,
    and it is a pinched set.
\item[\rm(3)] There is only one $\tau_\ep$-minimal set, in which case $\ml_\ep$ and $\muk_\ep$
    coincide on the residual set of common continuity points, and hence the global attractor
    is a pinched set. If this minimal set is hyperbolic, then it is an attractive hyperbolic
    $\tau_\ep$-copy of the base, given by $\{\ml_\ep\}=\{\muk_\ep\}$, and it coincides with $\mA_\ep$.
\end{itemize}
\end{teor}
\begin{proof}
The existence of the global attractor (which is compact and $\tau_\ep$-invariant) ensures
the existence of at least one $\tau_\ep$-minimal set: see Section \ref{subsec:minimal}.
According to \cite[Theorem 4.2]{duno1}, there are at most three of them, and the situation is
that of (1) if there are three. Let us define $\mM_\ep^l$ and $\mM_\ep^u$ from $\ml_\ep$ and
$\muk_\ep$ as in Proposition \ref{prop:defM}(i). Since any minimal set projects on the whole
set $\W$ (see Section \ref{subsec:minimal}), the existence of more than one of them ensures that
$\ml_\ep(\w)<\muk_\ep(\w)$ for all $\w\in\W$. So, if there are exactly two, then they are
$\mM_\ep^l$ and $\mM_\ep^u$: see again Proposition \ref{prop:defM}(i).
In addition, \cite[Theorem 5.13(iii)]{duno1} ensures that one of them is an attractive
hyperbolic copy of the base, and it follows from \cite[Proposition 5.3(ii)]{duno1}
that the other one is nonhyperbolic. So, we are in the situation (2). Finally, if there exists
exactly one minimal set, then $\mM_\ep^l=\mM_\ep^u$, and Proposition \ref{prop:defM}(i) guarantees
that $\ml_\ep(\w)=\muk_\ep(\w)$ at all the common continuity points of both maps; so, the section
$(\mA_\ep)_\w$ reduces to one element at these points,
and hence $\mA_\ep$ is pinched. If, in addition, $\mM_\ep^l=\mM_\ep^u$
is hyperbolic, then \cite[Proposition 5.3(i)]{duno1} precludes the possibility that it is repulsive, so
it is attractive. Hence, ``all" the minimal sets are hyperbolic attractive, which, according to
Theorem \ref{th:caract-hiper}(iii), ensures that $\mA_\ep$ is an attractive hyperbolic copy of the base.
That is, the situation is that described in (3).
\end{proof}
\begin{notas}\label{nota:Z}
1.~Whenever the dynamics of \eqref{eq:EDOpol}$_\ep$ fits in situation (1) of Theorem \ref{th:nummin},
we represent by $\{\mm_\ep\}$ the repulsive hyperbolic $\tau_\ep$-copy of the base.

2. A detailed description of the global dynamics (i.e., of the asymptotic behavior of the solutions)
can be done in each case of Theorem \ref{th:nummin}. We omit this, which is basically done
in \cite{duno1,duno4}, and we refer to the numerical simulations of Section \ref{sec:model}
for some clues in this regard.
\end{notas}
As said in the Introduction, we will perform our analysis in the case that $a(\w)<0$ and $c(\w)>0$ for all
$\w\in\W$. These are the unique conditions required in Theorem \ref{th:cpositivo} (and in the auxiliary
Proposition \ref{prop:raices}) to establish the first basic bifurcation properties.
\begin{prop}\label{prop:raices}
Assume that $a(\w)<0$ and $c(\w)>0$ for each $\w\in\W$. Then, there exists $\ep_0>0$ such that
\begin{itemize}[leftmargin=20pt]
\item[\rm(i)] the map $x\mapsto p_\ep(\w,x)$ has three real roots if $\ep\in(0,\ep_0]$ for all $\w\in\W$:
    $x^1_\ep(\w)>x^2_\ep(\w)>0>x^3_\ep(\w)$. In addition, $\lim_{\ep\to 0^+} (x^1_\ep(\w)-c(\w))=0$ and
    $\lim_{\ep\to 0^+}x^2_\ep(\w)=\lim_{\ep\to 0^+}x^3_\ep(\w)=0$, and the three limits are uniform on~$\W$.
\item[\rm(ii)] The map $x\mapsto p_\ep(\w,x)$ has only one real root if $\ep\in[-\ep_0,0)$, $x^1_\ep(\w)$, with
    $\lim_{\ep\to 0^-} (x^1_\ep(\w)-c(\w))=0$ uniformly on $\W$. In addition, $x^1_\ep(\w)>0$
    for all $\w\in\W$.
\end{itemize}
\end{prop}
\begin{proof}
A classical algebraic result (see, e.g., \cite[Exercises 10.14 and 10.17]{irvi}) establishes that the existence of one or three
real roots of the third degree polynomial $p_\ep(\w,x)$ depends on the sign of its {\em discriminant}
$\Delta_\ep(\w)$, given by
\begin{equation}\label{def:disc}
 \Delta_\ep:=\ep\,\big(-4\,a\,c^3+\ep\,b^2\,c^2 - 18\,\ep\,a\,b\,c-27\,\ep\,a^2+4\,\ep^2\,b^3\,\big):
\end{equation}
there is only one real root if $\Delta_\ep(\w)<0$ and three of them if $\Delta_\ep(\w)>0$. Hence,
since $\Delta_\ep(\w)$ is jointly continuous in $(\ep,\w)$, and since
$\lim_{\ep\to 0}\Delta_\ep(\w)/\ep=-4\,a(\w)\,c^3(\w)>0$, there are three real roots
$x^1_\ep(\w)>x^2_\ep(\w)>x^3_\ep(\w)$ if $\ep>0$ is small enough, and a real root
$x^1_\ep(\w)$ (plus two complex ones) if $-\ep>0$ is small enough. In addition,
the roots (considered as complex numbers) can be written as continuous maps
of the coefficients $\ep\,a$, $\ep\,b$ and $c$ of $p_\ep$.
Hence, the limits of the three solutions as $\ep\to 0$
are $c(\w)$, $0$ and $0$ (the roots of $x^2\,(c(\w)-x)=p_0(\w,x)=\lim_{\ep\to 0} p_\ep(\w,x)$),
and they are uniform on $\W$ (since the limits $\lim_{\ep\to 0}\ep\,a(\w)=\lim_{\ep\to 0}\ep\,b(\w)
=0$ are uniform on $\W$).
In both cases, the upper (unique for small $\ep<0$) real solution
converges to $c(\w)$ as $\ep\to 0$, so $x^1_\ep>0$ for $|\ep|>0$ small enough.
Finally, for $\ep>0$, $p_\ep(\w,0)=\ep\,a(\w)<0$ and $\lim_{x\to-\infty}p_\ep(\w,x)=\infty$.
So, if there were more than one negative root, there would be no positive ones, which as just
seen is precluded for $\ep>0$ small enough. The conclusion is $x^1_\ep(\w)>x^2_\ep(\w)>0>x^3_\ep(\w)$
for such an $\ep>0$.
\end{proof}
We fix some notation which will be used in the rest of the paper:
\[
\begin{split}
 &a_-:=\inf_{\w\in\W}a(\w) \qquad\text{and}\qquad a_+:=\sup_{\w\in\W}a(\w)\,,\\
 &b_-:=\inf_{\w\in\W}b(\w) \qquad\text{and}\qquad b_+:=\sup_{\w\in\W}b(\w)\,,\\
 &c_-:=\inf_{\w\in\W}c(\w) \qquad\text{and}\qquad c_+:=\sup_{\w\in\W}c(\w)\,.\\
\end{split}
\]
Recall that we say that there exists a (nonautonomous)
local saddle-node bifurcation point at $\ep_0$
when two hyperbolic copies of the base exist for $\ep<\ep_0$ (or $\ep>\ep_0)$ close to $\ep_0$ and they approach each
other as $\ep\to(\ep_0)^-$ (or as $\ep\to(\ep_0)^+$), giving rise to a locally unique nonhyperbolic $\tau_{\ep_0}$-minimal set $\mM_{\ep_0}$, and to the absence of minimal sets ``nearby $\mM_{\ep_0}\!$" for close $\ep>\ep_0$ (or $\ep<\ep_0)$.
\begin{teor}\label{th:cpositivo}
Assume that $a(\w)<0$ and $c(\w)>0$ for each $\w\in\W$. Then,
\begin{itemize}[leftmargin=20pt]
\item[\rm(i)] the unique $\tau_0$-minimal sets are $\{\ml_0\}=\{0\}$, which is nonhyperbolic,
    and $\{\muk_0\}$, which is hyperbolic attractive and satisfies $c_-\le \muk_0\le c_+$.
    In addition, either $\muk_0\equiv c$ and these maps are constant or, for all $\w\in\W$, there exists a
    strictly increasing two-sided sequence $(t_n)_{n\in\Z}$ with $\lim_{n\to\pm\infty}t_n=\pm\infty$
    such that $c(\wt_{2n})-\muk_0(\wt_{2n})>0$ and $c(\wt_{2n+1})-\muk_0(\wt_{2n+1})<0$.
\item[\rm(ii)] For all $\ep>0$, $\ml_\ep<0$, $\{\ml_\ep\}$ is an attractive hyperbolic $\tau_\ep$-copy of the base,
    and $(0,\infty)\to C(\W,\R),\,\ep\mapsto\ml_\ep$ is a continuous map in the uniform topology of $C(\W,\R)$.
\item[\rm(iii)] The set $\mI:=\{\ep_+>0\,|\;$there are three hyperbolic $\tau_\ep$-copies of the base
    for all $\ep\in(0,\ep_+)\}$ is nonempty and open; $\ml_\ep<0<\mm_\ep<\muk_\ep$ for all $\ep\in\mI$;
    $\lim_{\ep\to 0^+}(\muk_\ep(\w)-\muk_0(\w))=\lim_{\ep\to 0^+}\ml_\ep(\w)=\lim_{\ep\to 0^+}\mm_\ep(\w)=0$,
    all of them uniformly on $\W$;
    the maps $\mI\cup\{0\}\to C(\W,\R),\;\ep\mapsto\ml_\ep,\,\mm_\ep,\,\muk_\ep$ are continuous in the uniform
    topology of $C(\W,\R)$, where $\mm_0:=0$; and there exists
    $\ep_0\in(0,\sup\mI]$ such that the maps $(0,\ep_0)\to C(\W,\R),\,\ep\mapsto-\ml_\ep,\mm_\ep$
    are strictly increasing.
\item[\rm(iv)] For all $\ep<0$, $\muk_\ep>0$.
\item[\rm(v)] If, in addition, $c_+<3\,c_-$, then there exists
    $\ep_-<0$ such that, if $\ep\in(\ep_-,0)$,
    then $\mA_\ep=\{\muk_\ep\}$ is the unique $\tau_\ep$-minimal set; it is hyperbolic attractive;
    and $\lim_{\ep\to 0^-}(\muk_\ep(\w)-\muk_0(\w))=0$ uniformly on $\W$. In particular, there is a
    local saddle-node bifurcation at $\ep=0$.
\end{itemize}
\end{teor}
\begin{proof}
(i) Since $p_0(\w,x)=-x^3+c(\w)\,x^2$, the (unique)
Lyapunov exponent of the $\tau_0$-copy of the base $\{0\}$ is $\int_\W (p_0)_x(\w,0)\,dm=0$,
and hence $\{0\}$ is a nonhyperbolic $\tau_0$-minimal set: see Theorem \ref{th:caract-hiper}.
In particular, $\ml_0\le 0$. Since $p_0(\w,r)=-r^2(r-c(\w))>0$ for all $r<0$, Remark \ref{rm:comparison}.1
ensures that $\ml_0\ge0$, and hence $\ml_0\equiv 0$.
In addition, $p_0(\w,c_-)=(c_-)^2(c(\w)-c_-)\ge 0$ and
$p_0(\w,c_+)=(c_+)^2(c(\w)-c_+)\le 0$, and hence Proposition \ref{prop:encaje} ensures the existence
of a minimal set contained in $\W\times[c_-,c_+]\subset\W\times(0,\infty)$. So, we are necessarily
in case (2) of Theorem \ref{th:nummin}, and hence the second minimal set is $\{\muk_0\}$ and it is
hyperbolic attractive. The first assertions in (i) are proved.
\par
Now, observe that $\muk_0'(\wt)/\muk_0(\wt)=\muk_0(\wt)(-\muk_0(\wt)+c(\wt))$ for all $\w\in\W$ and $t\in\R$.
So, if $\muk_0\equiv c$, then $t\mapsto\muk_0(\wt)$ is constant for all $\w$ and hence
the continuous map $\w\mapsto\muk_0(\w)=c(\w)$ is constant by
the minimality of the base. Otherwise, Birkhoff's Ergodic Theorem yields
$0=\int_\W(\muk_0'(\w)/\muk_0(\w))\,dm=\int_\W \muk_0(\w)(-\muk_0(\w)+c(\w))\,dm=0$,
which precludes the global inequalities $c>\muk_0$ and $c<\muk_0$. Hence, the sets $\mU_\pm:=\{\w\in\W\,|\;
\pm(\muk_0(\w)-c(\w))>0\}$ are nonempty and open. The minimality of $\W$
yields, for a fixed $\w\in\W$, two increasing sequences
$(s_n^\pm)_{n\in\N}$ with limit $\infty$ and two decreasing ones $(\bar s_n^\pm)_{n\in\N}$
with limit $-\infty$ such that $\ws^\pm_n\in\mU_\pm$ and $\w{\cdot}\bar s^\pm_n\in\mU_\pm$,
which implies the existence of the sequence $(t_n)_{n\in\Z}$ of the last assertion in (i).
\smallskip\par
(ii)\&(iv) Since $p_\ep(\w,0)=\ep\,a(\w)$, we have $(0)'=0>p_\ep(\w,0)$ if $\ep>0$ and $(0)'=0<p_\ep(\w,0)$ if $\ep<0$
for all $\w\in\W$. So, Remark \ref{rm:comparison}.2 guarantees $\ml_\ep<0$ for $\ep>0$ (in (ii)) and property (iv).
To prove the second assertion in (ii), we fix $\ep>0$, define
\[
 q_\ep(\w,x):=\left\{\begin{array}{rl} -x^3+c(\w)\,x^2+\ep\,\big(b(\w)\,x+a(\w)\big)&\text{if $x\le 0$\,,}\\[.1cm]
 c(\w)\,x^2+\ep\,\big(b(\w)\,x+a(\w)\big)&\text{if $x>0$\,,}\end{array}\right.
\]
and consider the induced skew-product $\bar\tau_\ep(t,\w,x)=(\wt,\bar v_\ep(t,\w,x))$.
It is easy to check that $q_\ep$ is globally continuous and $C^2$ with respect to $x$,
and that its second derivative $(q_\ep)_{xx}$ is strictly positive. Hence,
$x\mapsto q_\ep(\w,x)$ is strictly convex for all $\w$, which ensures that
$y\mapsto -q_\ep(\w,y)$ is strictly concave for all $\w$. Let us consider the time-reversed flow
$\sigma^-\colon\R\times\W\to\W,\;(t,\w)\mapsto\w{\cdot}(-t)$. It is easy to check that $(\W,\sigma^-)$
is minimal and uniquely ergodic. Note that $-q_\ep(\w,0)=-\ep\,a(\w)>0$ for all $\w\in\W$,
and that there exists $r>0$ such that $-q_\ep(\w,\pm r)<0$ for all $\w\in\W$. Hence, Proposition \ref{prop:encaje}
applied to the skew-product flow $(\WR,\bar\tau_\ep^-)$ defined from $y'=-q_\ep(\w{\cdot}(-t),y)$ (over $(\W,\sigma^-)$)
ensures the existence of at least two minimal sets, strictly below and above $\W\times\{0\}$. Consequently, the dynamics for $\bar\tau_\ep^-$ is determined by an attractor-repeller pair of copies of the base (see,
e.g., \cite[Theorem 3.6]{duno4}), and this ensures that the lower minimal set is a repulsive hyperbolic $\bar\tau_\ep^-$-copy of the base, and that there are no bounded $\bar\tau_\ep^-$-orbits below it. The
change $x(t)=y(-t)$ takes $y'=-q_\ep(\w{\cdot}(-t),y)$ to $x'=q_\ep(\wt,x)$, and it is easy to deduce
from $\bar\tau_\ep^-(t,\w,x)=(\w{\cdot}(-t),\bar v_\ep(-t,\w,x))$
that the lower minimal set is an attractive hyperbolic $\bar\tau_\ep$-copy of the base, say $\{\bar\ml_\ep\}$,
and that there are no bounded $\bar\tau_\ep$-orbits below it.
Since each one of the maps $t\mapsto\bar\ml_\ep(\wt)$ and $t\mapsto\ml_\ep(\wt)$ (bounded and negative)
solve $x'=p_\ep(\wt,x)$ and $x'=q_\ep(\wt,x)$,
$\bar\ml_\ep\ge\ml_\ep$ (since $t\mapsto\ml_\ep(\wt)$ is the lower bounded solution for $x'=p_\ep(\wt,x)$), and
$\ml_\ep\ge\bar\ml_\ep$ (since $t\to\bar\ml_\ep(\wt)$ is the lower bounded solution for $x'=q_\ep(\wt,x)$).
Therefore, they are equal. In particular, $\int_{\W}p_x(\wt,\ml_\ep(\w))\,dm
=\int_{\W}q_x(\wt,\bar\ml_\ep(\w))\,dm$ (where $m$ is the unique $\sigma$-invariant measure on $\W$).
That is, the unique Lyapunov exponent for $\{\ml_\ep\}$ for the flow $\tau_\ep$ coincides
with that of $\{\bar\ml_\ep\}$ for $\bar\tau_\ep$, and hence Theorem \ref{th:caract-hiper}
ensures that it is negative (since it is the unique Lyapunov exponent of the attractive hyperbolic
$\bar\tau_\ep$-copy of the base $\{\bar\ml_\ep\}$) and that $\{\ml_\ep\}$ is an attractive
hyperbolic $\tau_\ep$-copy of the base (since its unique Lyapunov exponent is negative).
\par
Finally, the classical result of robustness of the existence of hyperbolic copies of the base
and their continuous variation in the uniform topology (see, e.g., \cite[Theorem 2.3]{duno5})
proves the continuity of $\ep\mapsto\ml_\ep$ stated in (ii).
\smallskip\par
(iii) Let us take $\rho\in(0,c_-)$. Proposition \ref{prop:raices}(i) allows us to take $\ep_\rho>0$
small enough to ensure that $x^3_\ep(\w)<0<x^2_\ep(\w)<\rho<x^1_\ep(\w)$
for all $\w\in\W$ and $\ep\in(0,\ep_\rho)$.
Since $p_\ep(\w,x)=-(x-x^1_\ep(\w))(x-x^2_\ep(\w))(x-x^3_\ep(\w))$, we have $p_\ep(\w,\rho)>0$ for all $\w\in\W$.
We also look for $k_-<0$ and $k_+>\rho$ such that $p_\ep(\w,k_-)>0$ and $p_\ep(\w,k_+)<0$ for
all $\ep\in[0,\ep_\rho]$ and all $\w\in\W$. Then: $k_-$ and $\rho$ are strict global upper solutions, and
$0$ and $k_+$ are strict global lower solutions. Since $k_-<0<\rho<k_+$, Proposition \ref{prop:encaje}
ensures the existence of three minimal sets $\mM^1_\ep$, $\mM^2_\ep$ and $\mM^3_\ep$,
with $\mM^1_\ep\subset\W\times(k_-,0)$, $\mM^2_\ep\subset\W\times(0,\rho)$ and
$\mM^3_\ep\subset\W\times (\rho,k_+)$. Theorem \ref{th:nummin} ensures that
$\mM^1_\ep=\{\ml_\ep\}$ (attractive), $\mM^2_\ep=\{\mm_\ep\}$ (repulsive), and
$\mM^3_\ep=\{\muk_\ep\}$ (attractive). So, $(0,\ep_\rho)\subseteq\mI$, with
$\mI$ defined in (iii). The robustness of the existence of hyperbolic copies of the base
and their continuous variation in the uniform topology (see \cite[Theorem 2.3]{duno5}) prove that
$\mI$ is open and that the maps $\mI\to C(\W,\R),\;\ep\mapsto\ml_\ep,\,\mm_\ep,\,\muk_\ep$
are continuous. Note also that, if $\ep>0$ the inequality $p_\ep(\w,0)<0$ for all $\w$
precludes the existence of a point $(\w,0)$ in any $\tau_\ep$-minimal set $\mM_\ep$:
see Remark \ref{rm:comparison}.3. So, since $\mm_\ep>0$ for small $\ep>0$ and
$\ep\to\mm_\ep$ is continuous on $\mI$, we conclude that $\mm_\ep>0$ and hence $\ml_\ep<0<\mm_\ep<\muk_\ep$
for all $\ep\in\mI$.
\par
Keeping the notation of Proposition \ref{prop:raices}, we define $(x^3_\ep)_-:=\inf_{\w\in\W}x^3_\ep(\w)$
for $\ep>0$ small enough. Then, $p_\ep(\w,r)>0$ for all $r<(x^3_\ep)_-$, which ensures that
$\ml_\ep>(x^3_\ep)_-$: see Remark \ref{rm:comparison}.1.
Proposition \ref{prop:raices} ensures that $\lim_{\ep\to 0^+}(x^3_\ep)_-=0$, which combined with $l_\ep<0$
for $\ep>0$ yields $\lim_{\ep\to 0^+}\ml_\ep(\w)=0$ uniformly on $\W$. Now, we choose the initial $\rho>0$
as close to $0$ as desired, and observe that $0<\mm_\ep<\rho$ if $\ep$ is small enough. This proves the assertion
concerning $\lim_{\ep\to 0^+}\mm_\ep$. Finally, \cite[Theorem 2.3]{duno5} ensures the existence of an
attractive hyperbolic copy of the base as uniformly close as desired to $\{\muk_0\}$ for small $\ep>0$,
which must be $\{\muk_\ep\}$, and this proves the assertion about $\lim_{\ep\to 0^+}(\muk_\ep-\muk_0)$.
These properties and (i) complete the proof of the continuity on $\mI\cup\{0\}$.
\par
The proofs of the monotonicity properties of $\ep\mapsto\ml_\ep,\mm_\ep$
take arguments from the proof of \cite[Theorem 5.10]{duno1}, which
we repeat here for the reader's convenience. The previous uniform limiting properties allow us to ensure that,
if $\ep>0$ is small enough, then $b(\w)\,x+a(\w)<0$ if $x\in[\ml_\ep(\w),\mm_\ep(\w)]$ for all $\w\in\W$.
We choose $\ep_0$ so that this holds for $\ep\in(0,\ep_0)$.
If $0<\ep_1<\ep_2<\ep_0$, then $\ml'_{\ep_1}(\w)=p_{\ep_1}(\w,\ml_{\ep_1}(\w))>
p_{\ep_2}(\w,\ml_{\ep_1}(\w))$, and hence $\ml_{\ep_1}>\ml_{\ep_2}$ (see Remark \ref{rm:comparison}.2).
Let us complete the proof of (iii) checking that $\ep\mapsto\mm_\ep$ is strictly increasing on $(0,\ep_0)$.
We fix $\ep_1\in(0,\ep_0)$. The previously checked continuous variation of the copies of the base
allows us to take $\ep_2<\ep_1$ in $(0,\ep_0)$ close enough to ensure $\ml_{\ep_1}<\mm_{\ep_2}<\muk_{\ep_1}$.
So, for a fixed $\bar\w\in\W$, $v_{\ep_1}(t,\bar\w,\mm_{\ep_2}(\bar\w))>v_{\ep_1}(t,\bar\w,\ml_{\ep_1}(\bar\w))=
\ml_{\ep_1}(\bwt)$ for all $t\in\R$: the $\tau_{\ep_1}$-orbit of $(\bar\w,\mm_{\ep_2}(\bar\w))$ is
above $\{\ml_{\ep_1}\}$, and hence globally bounded. This ensures the existence of the
corresponding \upalfa-limit set. In addition, $\mm_{\ep_2}'(\w)<p_{\ep_1}(\w,\mm_{\ep_2}(\w))$ for all $\w\in\W$, since
${\ep_2}\,(b(\w)\,\mm_{\ep_2}(\w)+a(\w))<{\ep_1}\,(b(\w)\,\mm_{\ep_2}(\w)+a(\w))$.
According to Remark \ref{rm:comparison}.3, the previous \upalfa-limit set is strictly below
the graph of $\mm_{\ep_2}$.
Let $\mN$ be a $\tau_{\ep_1}$-minimal contained in this \upalfa-limit set.
This \upalfa-limit set cannot intersect $\{\ml_{\ep_1}\}$ or $\{\muk_{\ep_1}\}$
(see Remark \ref{rm:comp2}), and hence $\mN=\{\mm_{\ep_1}\}$. This ensures that
$\mm_{\ep_1}<\mm_{\ep_2}$, as asserted.
\smallskip\par
(v) We assume that $c_+<3\,c_-$ and work with $\ep<0$. Initially, we take values of $\ep<0$ close enough to $0$
to ensure that, for all $\w\in\W$, $x^1_\ep(\w)>0$ is the unique real root of $x\mapsto p_\ep(\w,x)$
(see Proposition \ref{prop:raices}), and define $(x_\ep^1)_-:=\inf_{\w\in\W}x^1_\ep(\w)$.
Then $\ml_\ep>(x_\ep^1)_-$, since $p_\ep(\w,r)>0$ for all $r<(x_\ep^1)_-$ (see again
Remark \ref{rm:comparison}.1). Since $\lim_{\ep\to 0^-} (x_\ep^1(\w)-c(\w))=0$ uniformly on $\W$
(see again Proposition \ref{prop:raices}), we have $\lim_{\ep\to 0^-}(x_\ep^1)_-=c_-$.
We fix $r\in(c_+/3,c_-)$ and look for $\ep_-<0$ close enough to $0$ to ensure that
$(x^1_\ep)_->r$ for all $\ep\in(\ep_-,0)$. Now, we fix $\ep\in(\ep_-,0)$ and define
$q_\ep(\w,x)$ as the $C^{0,2}(\WR,\R)$ function which coincides with $p_\ep(\w,x)$ for $x\ge r$
and is given by a second degree polynomial for $x\le r$. In particular, $q_\ep(\w,r)>0$. In addition,
for all $x\le r$, $(\partial^2/\partial x^2)\,q_\ep(\w,x)=(\partial^2/\partial x^2)\,q_\ep(\w,r)=
(\partial^2/\partial x^2)\,p_\ep(\w,r)=-6\,r+2\,c(\w)\le-6\,r+2\,c_+<0$. And, if $x>r$, then
$(\partial^2/\partial x^2)\,q_\ep(\w,x)=(\partial^2/\partial x^2)\,p_\ep(\w,x)=
-6\,x+2\,c(\w)<-6\,r+2\,c_+<0$. That is, the map $x\mapsto q_\ep(\w,x)$ is strictly concave
for all $\w\in\W$. Moreover, $\lim_{x\to-\infty}q_\ep(\w,x)=-\infty$ (as in the case of a concave
second-degree polynomial) and $\lim_{x\to\infty}q_\ep(\w,x)=\lim_{x\to\infty}p_\ep(\w,x)=-\infty$,
uniformly on $\W$ in both cases. This properties mean that $q_\ep$ satisfies all the conditions
{\bf c1}-{\bf c4} of \cite[Section 3]{duno4}. Since $\ml_\ep\ge (x^1_\ep)_->r$,
any $\tau_\ep$-minimal set is contained in $\W\times(r,\infty)$,
and hence is also a $\tilde\tau_\ep$-minimal set for
the skew-product flow $\tilde\tau_\ep$ defined on $\W\times\R$ by $x'=q_\ep(\wt,x)$. Conversely,
any $\tilde\tau_\ep$-minimal set contained in $\W\times(r,\infty)$ is also a $\tau_\ep$-minimal set.
Since $q_\ep(\w,r)>0$
and $q_\ep(\w,\pm s)<0$ for all $\w\in\W$ if $s$ is large enough, Proposition \ref{prop:encaje} ensures
the existence of a $\tilde\tau_\ep$-minimal set $\mM^u_\ep$ in $\W\times(r,\infty)$
and of another in $\W\times(-\infty,r)$.
According to \cite[Theorem 3.3]{duno4}, $\mM^u_\ep$ is an attractive hyperbolic $\tilde\tau_\ep$-copy
of the base, and the unique $\tilde\tau_\ep$-minimal set above $r$. The conclusion is that
$\mM^u_ \ep$ is also the unique $\tau_\ep$-minimal set. Since its Lyapunov exponents are the same
for $p_\ep$ as for $q_\ep$ (due to the equality $(p_\ep)_x=(q_\ep)_x$ on $\W\times(r,\infty)$),
they are negative, and Theorem \ref{th:caract-hiper} ensures that $\mA_\ep=\mM_\ep$ is
an attractive hyperbolic $\tau_\ep$-copy of the base, as asserted. Finally,
the last assertion in (v) follows from the previous ones and (iii).
\end{proof}
\begin{notas}\label{notas:ZZ}
1.~Observe that the hypothesis $c_+<3\,c_-$ of Theorem \ref{th:cpositivo}(v) always
holds if the map $c$ is a positive constant. Otherwise, the range of
``allowed" values of $c$ increases as $c_-$ increases.
Note also that, if all the conditions assumed in Theorem \ref{th:cpositivo}(v) hold, then
the local saddle-node bifurcation at $\ep=0$ has an extra property:
the collision as $\ep\to 0^+$ of the two approaching hyperbolic
$\tau_\ep$-copies of the base is total, giving rise to the nonhyperbolic $\tau_0$-copy of the base $\{0\}$.
So, in contrast to the possibly  very complex dynamics at the bifurcation values $\ep\ne 0$ which we will find later
(due to the possibly very complex dynamics of the nonhyperbolic minimal set),
here the dynamics for $\ep=0$ is simple: a nonautonomous reproduction of the autonomous dynamics around $\ep=0$
of, for instance, $x'=-x^3+x^2-\ep$. In fact, since $p_0(\w,r)>0$ for all $r\in(0,c_-)$, it is easy to check that
$\{0\}$ is the \upalfa-limit set of the $\tau_0$-orbit of all $(\w,r)\in\W\times[0,c_-)$ and that $\{\muk_0\}$ is the unique
$\tau_0$-minimal set contained in the \upomeg-limit set.
\par
2. If, under the hypotheses of Theorem \ref{th:cpositivo}, we also assume $b\equiv 0$,
then an easy extension of the results of \cite{duno1}
(see also \cite{duno5}) provides a complete description of the global bifurcation diagram of \eqref{eq:EDOpol}:
\begin{itemize}[leftmargin=10pt]
\item[-] there are exactly two bifurcation points, $0$ and a certain $\ep_*>0$;
\item[-] the maps $(-\infty, \ep_*)\to C(\W,\R),\,\ep\mapsto\muk_\ep$ and $(0,\infty)\to C(\W,\R),\,\ep\mapsto\ml_\ep$
are continuous and strictly decreasing, and $\{\muk_\ep\}$ (resp.~$\{\ml_\ep\}$)
is an attractive hyperbolic $\tau_\ep$-copy of the base for all $\ep<\ep_*$ (resp. $\ep>0$);
\item[-] $\mA_\ep=\{\ml_\ep\}=\{\muk_\ep\}$ for $\ep\notin[0,\ep_*]$, with
$\lim_{\ep\to\pm\infty}\muk_\ep=\mp\infty$ uniformly on $\W$;
\item[-] there are three hyperbolic $\tau_\ep$-copies of the base
for $\ep\in(0,\ep_*)$ given by $\ml_\ep<\mm_\ep<\muk_\ep$, and $\ep\mapsto\mm_\ep$
is strictly increasing on $(0,\ep_*)$;
\item[-] and there are two $\tau_\ep$ minimal sets for $\ep=0$ (resp.~$\ep=\ep_*$):
$\{\muk_0\}$ (resp.~$\{\ml_{\ep_*}\}$), which is hyperbolic attractive,
and a nonhyperbolic one given by the collision of the two lower (resp.~upper)
copies of the base as $\ep\to 0^+$ (resp. $\ep\to(\ep_*)^-$).
\end{itemize}
So, the two bifurcations points are of local saddle-node type. The interested reader can find in
\cite[Figure 1]{duno1} a similar bifurcation diagram, which must be horizontally inverted to get ours.
In this situation, Theorem \ref{th:cpositivo}
adds just a little piece of information to these facts: as explained in the previous remark,
the lower nonhyperbolic minimal set for $\ep=0$ reduces to the copy of the base $\{0\}$.
Note also that this description shows that the situation that Theorem \ref{th:cbpositivo}(iii) describes
cannot hold if $b\equiv 0$.
\end{notas}
The hypotheses $b\ge 0$ or $b>0$, added in the next result, allow us to delve deeper into
the dynamical changes as $\ep$ varies.
\begin{teor}\label{th:cbpositivo}
Assume that $a(\w)<0$, $b(\w)\ge 0$ and $c(\w)>0$ for all $\w\in\W$. Then,
in addition to the information provided by Theorem {\rm\ref{th:cpositivo}},
\begin{itemize}[leftmargin=20pt]
\item[\rm(i)] $\ml_\ep<0$ for all $\ep>0$, $\ml_0\equiv 0$, and $\ml_\ep>s_0>0$ for all $\ep<0$,
where $s_0:=\min(\,\inf_{\w\in\W}c(\w)\,,\,\inf_{\w\in\W,\,b(\w)\ne 0}(-a(\w)/b(\w))\,)$ if $b\not\equiv 0$, and
$s_0:=\inf_{\w\in\W}c(\w)$ if $b\equiv 0$. In particular, there is a local saddle-node bifurcation at $\ep=0$.
\item[\rm(ii)] The continuous map $(0,\infty)\to C(\W,\R),\,\ep\mapsto\ml_\ep$ is strictly decreasing,
    with $\lim_{\ep\to 0^+}\ml_\ep(\w)=0$ and $\lim_{\ep\to \infty}\ml_\ep(\w)=-\infty$, both of them uniformly on $\W$.
\end{itemize}
If, in addition, $b(\w)>0$ for all $\w\in\W$, then
\begin{itemize}[leftmargin=20pt]
\item[\rm(iii)] the set $\mJ:=\{\ep^+>0\,|\;$ there are three hyperbolic $\tau_\ep$-copies of the base for all $\ep>\ep^+\}$
is nonempty and open, the maps $\mJ\to C(\W,\R),\;\ep\mapsto\ml_\ep,\,\mm_\ep,\,\muk_\ep$
are continuous in the uniform topology of $C(\W,\R)$, and $\lim_{\ep\to \infty}\muk_\ep(\w)=\infty$
uniformly on $\W$. In addition, there exists $\ep^0\ge\inf\mJ$ such that
the continuous map $(\ep^0,\infty)\to C(\W,\R),\,\ep\mapsto\muk_\ep$ is strictly increasing.
\end{itemize}
\end{teor}
\begin{proof}
(i) The two first assertions in (i) (as well as the continuity stated in (ii))
are proved in Theorem \ref{th:cpositivo}(i)-(ii).
Let us now define $s_0$ as in the statement of (i), and take $r<s_0$. Then, if $\ep<0$,
$p_{\ep}(\w,r)>p_0(\w,r)>0$: the first inequality follows from $\ep\,(b(\w)\,r+a(\w))>0$, and the second one
from $c(\w)-r>0$. This ensures that
$\mA_\ep\subset\W\times(s_0,\infty)$ (see Remark \ref{rm:comparison}.1) and hence
proves that $\ml_\ep>s_0$ for all $\ep<0$.
This fact precludes the existence of $\tau_\ep$-minimal sets ``close" to $\{0\}$ if $\ep<0$. Consequently,
there is a local saddle-node bifurcation at $\ep=0$, due to the collision of $\ml_\ep$ and $\mm_\ep$
as $\ep\to 0^+$: see Theorem \ref{th:cpositivo}(iii).
\smallskip\par
(ii) Since $\ml_\ep<0$ for $\ep>0$, we have $b(\w)\,\ml_{\ep}(\w)+a(\w)<0$ for all $\w\in\W$. So, if $0<\ep_1<\ep_2$, then
$\ml_{\ep_1}'(\w)=p_{\ep_1}(\w,\ml_{\ep_1}(\w))>p_{\ep_2}(\w,\ml_{\ep_1}(\w))$,
and hence $\ml_{\ep_1}>\ml_{\ep_2}$: see Remark \ref{rm:comparison}.2.
Theorem \ref{th:cpositivo}(v) proves the assertion about $\lim_{\ep\to 0^+}\ml_\ep$.
Now, we fix $r<0$
and note that $\lim_{\ep\to\infty} p_\ep(\w,r)=-\infty$ uniformly on $\W$, since $a(\w)+b(\w)\,r\le a_+<0$. Hence,
there exists $\ep_r>0$ such that $p_\ep(\w,r)<0$ for all $\ep\ge\ep_r$. According to Remark \ref{rm:comparison}.2,
$\ml_\ep<r$ for all $\ep>\ep_r$, which proves the last assertion in~(ii).
\smallskip\par
(iii) The first goal is to find a strictly positive constant providing a global strict lower solution
if $\ep>0$ is large enough: as we will see later, the conclusions follow from this property.
Note that, if $\ep>0$ and $x>0$, then $p_\ep(\w,x)\ge\bar p_\ep(x):=-x^3+c_-\,x^2+\ep\,b_-\,x+\ep\,a_-$ for all $\w\in\W$.
The points $x^+_\ep$ and $x^-_\ep$ given by
\[
 x^\pm_\ep=\frac{c_-\pm\sqrt{c_-^2+3\,\ep\,b_-}}3
\]
are the local minimum and maximum of $\bar p_\ep(x)$, respectively.
We observe that $x^+_\ep>0$ and write it as
$x^+_\ep=(c_-+r_\ep)/3$, with $r_\ep:=\sqrt{c_-^2+3\,\ep\,b_-}>\sqrt{3\,\ep\,b_-}$.
A straightforward computation shows that
\[
\begin{split}
 27\,\bar p_\ep(x^+_\ep)&=6\,\ep\,b_-\,r_\ep+9\,\ep\,b_-\,c_-+27\,\ep\,a_-+2\,c_-^2\,r_\ep+2\,c_-^3\\
 &>6\sqrt{3}\,\ep^{3/2}\,b_-^{3/2} + \ep\,(27\,a_-+9\,b_-\,c_-)\,,
\end{split}
\]
and hence $p_\ep(\w,x^+_\ep)\ge\bar p_\ep(x^+_\ep)>0$ for all $\w\in\W$ if $\ep>0$ is large enough.
\par
So, for $\ep>0$ large enough, the constant $x^+_\ep$ is a global strict lower solution.
The expression of $p_\ep(\w,x)$ shows that a sufficiently large constant $s_\ep>0$ is a global
strict upper solution, and $0$ is a global strict upper solution (due to $p_\ep(\w,0)=\ep\,a(\w)<0$).
So, Proposition \ref{prop:encaje} ensures the existence of two minimal sets:
$\mM_\ep^u\subset\W\times(x^+_\ep,\infty)$ and $\mM_\ep^m\subset\W\times(0,x^+_\ep)$.
Since Theorem \ref{th:cpositivo}(ii) shows the existence of a third minimal set $\mM_\ep^l=\{\ml_\ep\}
\subset\W\times(-\infty,0)$, Theorem \ref{th:nummin} ensures that $\mM_\ep^u=\{\muk_\ep\}$ and it is hyperbolic attractive,
and that $\mM_\ep^m$ is a repulsive hyperbolic $\tau_\ep$-copy of the base, say $\mM_\ep^m=\{\mm_\ep\}$.
\par
Therefore, $\mJ$ is nonempty. It is open, as a consequence
of the persistence of hyperbolic copies of the base under small variations of $\ep$
(see again \cite[Theorem 2.3]{duno5}), which also shows the continuity asserted in (iii). In addition,
$\lim_{\ep\to\infty}\muk_\ep(\w)=+\infty$ uniformly on $\W$, since $\muk_\ep>x^+_\ep$ and
$\lim_{\ep\to\infty}x^+_\ep=\infty$. This last property ensures the existence of $\ep^0\ge\inf\mJ$
such that $\inf_{\w\in\W}\muk_\ep(\w)\ge\sup_{\w\in\W}(-a(\w)/b(\w))$
for $\ep>\ep^0$. So, if $\ep_2>\ep_1>\ep^0$, then $\muk_{\ep_1}'(\w)=p_{\ep_1}(\w,\muk_{\ep_1}(\w))<
p_{\ep_2}(\w,\muk_{\ep_1}(\w))$ for all $\w\in\W$, and hence $\muk_{\ep_1}<\muk_{\ep_2}$ (see
Remark \ref{rm:comparison}.2). The proof is complete.
\end{proof}
\subsection{Four different nonautonomous bifurcation diagrams}
The next three results add extra conditions to $a<0$, $b\ge 0$ with $b\not\equiv 0$ and $c>0$, related
to the relative sizes of $b\,c\ge 0$ and $a<0$.
These conditions can be considered nonautonomous extensions of the three
possibilities arising in the autonomous case, namely
$b\,c+a<0$, $b\,c+a>0$ and $b\,c+a=0$.
The three cases which we consider are far away from covering the infinitely many
possibilities that arise in the nonautonomous case, but there are mathematical models
that justify their interest, such as the case of the population dynamics that we analyze in
Section \ref{sec:model}.
Recall that we assume the existence of a unique $\sigma$-invariant measure $m$ on $\W$.
\begin{teor}\label{th:caso1}
Assume that $a(\w)<0$, $b(\w)\ge 0$, $c(\w)>0$, and $b(\w)\,c_++a(\w)<0$ for all $\w\in\W$.
Let $\mI=(0,\ep_*)$ be the open interval of Theorem {\rm\ref{th:cpositivo}(iii)}. Then, in
addition to the information provided by Theorems {\rm\ref{th:cpositivo}} and {\rm\ref{th:cbpositivo}},
\begin{itemize}[leftmargin=20pt]
\item[\rm(i)] $\ep_*\in\R$, the continuous maps
    $\mI\to C(\W,\R),\,\ep\mapsto \ml_\ep,\,-\mm_\ep,\,\muk_\ep$ are strictly decreasing, and
    $\mm_{\ep_*}(\w):=\lim_{\ep\to(\ep_*)^-}\mm_\ep(\w)$ and
    $\bar\muk_{\ep_*}(\w):=\lim_{\ep\to(\ep_*)^-}\muk_\ep(\w)$ define two semicontinuous $\tau_{\ep_*}$-equilibria
    which coincide with $\muk_{\ep_*}(\w)$ at all $\w$ in a $\sigma$-invariant residual subset
    $\mR_{\ep_*}\!\subseteq\W$. In particular, there exist exactly two $\tau_{\ep_*}$-minimal sets:
    $\{l_{\ep_*}\}\subset\W\times(-\infty,0)$, which is hyperbolic attractive,
    and $\mM_{\ep_*}\subset\W\times(0,c_+)$, which is nonhyperbolic.
\item[\rm(ii)] For all $\ep<0$, $\mA_\ep\subset\W\times(c_-,\infty)$.
\item[\rm(iii)] If, in addition, $c_+<3\,c_-$, then $\mA_\ep=\{\muk_\ep\}$ is the unique $\tau_\ep$-minimal set for all $\ep<0$,
    it is hyperbolic attractive,
    and the map $(-\infty,\ep_*)\to C(\W,\R),\,\ep\mapsto\muk_\ep$ is continuous in the uniform topology of $C(\W,\R)$.
\end{itemize}
Assume also that $b(\w)> 0$ for all $\w\in\W$, and define $s_-:=\inf_{\w\in\W}(-a(\w)/b(\w))$ and
$s_+:=\sup_{\w\in\W}(-a(\w)/b(\w))$. Then,
\begin{itemize}[leftmargin=20pt]
\item[\rm(iv)] for all $\ep\le 0$, $\muk_\ep<s_+$; there exists $\bar\ep\le 0$ such that
    $\mA_\ep=\{\muk_\ep\}$ is an attractive hyperbolic copy of the base for all $\ep<\bar\ep$;
    and if, in addition, $c_+<3\,c_-$, then there exists $\underline\ep$ with
    $-\infty\le\underline\ep<0$ such that $\ep\mapsto\muk_\ep$ is strictly decreasing on $(\underline\ep,0)$.
\item[\rm(v)] Let $\mJ=(\ep^*,\infty)$ be the open set defined in
    Theorem {\rm\ref{th:cbpositivo}(iii)}. Then, $\ep_*<\ep^*$, $\mm_\ep>s_-$
    for all $\ep\in\mJ$, and there are exactly two $\tau_{\ep^*}$-minimal sets:
    $\{\ml_{\ep^*}\}$ and $\mM_{\ep^*}\subset\W\times(s_-,\infty)\subset\W\times(c_+,\infty)$.
\item[\rm(vi)] If, in addition, $a/b=-s\in\R$, then the maps $(\ep^*,\infty)\to C(\W,\R),\,
    \ep\mapsto\muk_\ep,\,-\mm_\ep$ are strictly increasing; the maps $\mm_{\ep^*}(\w):=\lim_{\ep\to(\ep^*)^+}\mm_\ep(\w)$ and
    $\bar\muk_{\ep^*}(\w):=\lim_{\ep\to(\ep^*)^+}\muk_\ep(\w)$ define two
    semicontinuous $\tau_{\ep^*}$-equilibria which coincide with $\muk_{\ep^*}(\w)$
    at all $\w$ in a $\sigma$-invariant residual subset $\mR_{\ep^*}\!\subseteq\W$;
    the continuous map $(-\infty,\ep_*)\to C(\W,\R),\,\ep\mapsto\muk_\ep$ is strictly decreasing;
    $\lim_{\ep\to-\infty}\muk_\ep(\w)=\lim_{\ep\to\infty}\mm_\ep(\w)=s$ uniformly on $\W$;
    and $\{\ml_\ep\}$ is the unique $\tau_\ep$ minimal set for $\ep\in(\ep_*,\ep^*)$.
    So, $0,$ $\ep_*$ and $\ep^*$ are three bifurcation points,
    all of them of local saddle-node type, and they are the unique ones if $c_+<3\,c_-$.
\end{itemize}
\end{teor}
\begin{proof}
(i) In what follows, we use the notation and information of Proposition \ref{prop:raices}.
Let us take $\ep_0\in\mI$ small enough to ensure that, if $\ep\in(0,\ep_0]$, then there exist three
real roots of $p_\ep(\w,x)$ which satisfy $x^3_\ep(\w)<0<x^2_\ep(\w)<x^1_\ep(\w)$ and $x^2_\ep(\w)<c_+$
for all $\w\in\W$, and take $\ep\in(0,\ep_0]$. Since
$p_\ep(\w,c_+)=(c_+)^2(c(\w)-c_+)+\ep\,(b(\w)\,c_++a(\w))<0$ if $\w\in\W$, either
$c_+<x^2_\ep(\w)$ or $c_+>x^1_\ep(\w)$ for each $\w\in\W$, so that $c_+>x^1_\ep(\w)$ for all $\w\in\W$.
So, $p_\ep(\w,r)<0$ for any $r>c_+$, and hence $\muk_\ep<c_+$: see Remark
\ref{rm:comparison}.1. Now assume that $\muk_\ep\le c_+$ for $\ep=\ep_1,\ep_2$ with $0\le\ep_1<\ep_2$.
Then, $\muk_{\ep_2}'(\w)=p_{\ep_2}(\w,\muk_{\ep_2}(\w))>p_{\ep_1}(\w,\muk_{\ep_2}(\w))$,
since $b(\w)\,\muk_\ep(\w)+a(\w)\le b(\w)\,c_++a(\w)<0$. Hence, $\muk_{\ep_2}<\muk_{\ep_1}\le c_+$
(see Remark \ref{rm:comparison}.2). In particular, $\ep\mapsto\muk_\ep$ strictly decreases on $[0,\ep_0]$.
Let $\mI_u\subset\mI$ the interval of persistence of this property. There exists $\delta>0$
such that, for any $\ep\in\mI_u$, $\muk_\ep<\muk_0-\delta\le c_+-\delta$. This, the continuity of
$\ep\mapsto\muk_\ep$ on $\mI$, and the previously proved property, preclude the possibility
that $\sup\mI_u<\sup\mI$, and hence $\ep\mapsto\muk_\ep$ strictly decreases on $\mI$.
\par
The strictly decreasing character of the continuous map $\ep\mapsto\ml_\ep$ on $\mI\subseteq(0,\infty)$
is proved in Theorem \ref{th:cbpositivo}(ii). To check that $\ep\mapsto\mm_\ep$ is strictly increasing on
$\mI$, we adapt the argument of the proof of Theorem \ref{th:cpositivo}(iii). First note that,
for all $\ep\in\mI$, $\mm_\ep<\muk_\ep<c_+$ and hence $b(\w)\,\mm_\ep(\w)+a(\w)<0$
for all $\w\in\W$.  So: we fix $\ep_1\in\mI$ and take
${\ep_2}>{\ep_1}$ in $\mI$ close enough to ensure $\muk_{\ep_1}>\mm_{\ep_2}>\ml_{\ep_1}$, so that
the $\tau_{\ep_1}$-orbit of $(\w,\mm_{\ep_2}(\w))$ is above $\{\ml_{\ep_1}\}$; we check that
$\mm_{\ep_2}$ is a global strict lower solution for $\tau_{\ep_1}$; we use these properties and
Remark \ref{rm:comparison}.3 to ensure that the
\upalfa-limit set for $\tau_{\ep_1}$ of $(\w,\mm_{\ep_2}(\w))$ exists
and is strictly below the graph of $\mm_{\ep_2}$; and we deduce that the unique $\tau_{\ep_1}$-minimal set contained in
this \upalfa-limit set is $\{\mm_{\ep_1}\}$, so that $\mm_{\ep_1}<\mm_{\ep_2}$.
\par
Now, we assume for contradiction that $\sup\mI=:\ep_*=\infty$. Note that $\muk_\ep>\mm_\ep>0$ for all
$\ep\in\mI$. Given any $\rho>0$, we take
$\ep_\rho>\max_{\w\in\W,\,x\in[-\rho,c_+]}(x^3-c(\w)\,x^2)/(b(\w)\,x+a(\w))\ge 0$
(note that $b(\w)\,x+a(\w)<0$ if $x\le c_+$).
Then, $p_\ep(\w,r)<0$ for all $r\in[-\rho,c_+]$ if $\ep>\ep_\rho$. Let us deduce
that $\muk_{\ep}\le-\rho<0$ if $\ep\ge\ep_\rho$, which provides the sought-for contradiction.
Recall that $\muk_\ep$ is continuous and $\muk_\ep<c_+$ for all $\ep\in\mI=(0,\infty)$.
Again for contradiction, we assume that $\max_{\w\in\W}\muk_\ep(\w)=:s\in(-\rho,c_+)$,
and take $\bar\w\in\W$ with $\muk_\ep(\bar\w)=s$. Then, $p_\ep(\bar\w,\muk_\ep(\bar\w))<0$,
and hence $\muk_\ep(\bar\w{\cdot}t)=v_\ep(t,\bar\w,\muk_\ep(\bar\w))>s$ for small values
of $t<0$, which contradicts the definition of $s$. The conclusion is
that $\ep_*$ is finite, as asserted.
\par
Recall that $\ep_*\notin\mI$, since $\mI$ is open. Therefore, there exists at most a $\tau_{\ep_*}$-minimal
set $\mM_{\ep_*}$ different from $\{\ml_{\ep_*}\}$ (see Theorems \ref{th:nummin} and \ref{th:cpositivo}(ii)),
which is nonhyperbolic. The remaining part of this proof is very similar to part of that of
\cite[Theorem 5.10(i)]{duno1}, but we detail it for the reader's convenience.
The monotonicity properties of $0<\mm_\ep<\muk_\ep$ ensure the global existence of the limits
$\mm_{\ep_*}$ and $\bar\muk_{\ep_*}$, with $0<\mm_\ep<\mm_{\ep_*}\le\bar\muk_{\ep_*}<\muk_\ep$ for $\ep\in\mI$.
It is easy to check that they are $\tau_{\ep_*}$-equilibria. Since they are monotone limits of continuous functions,
they are semicontinuous on $\W$: $\mm_{\ep_*}$ is lower semicontinuous and $\bar\muk_{\ep_*}$
is upper semicontinuous.
Let $\mM_{\ep_*}^m$ be the minimal set associated to $\mm_{\ep_*}$ by \eqref{eq:precollision}.
The lower semicontinuity of $\mm_{\ep_*}$ yields $x\ge\mm_{\ep_*}(\w)\ge 0$ for any $(\w,x)\in\mM$.
In particular, $\mM_{\ep_*}^m\neq\{\ml_{\ep_*}\}$ (see Theorem \ref{th:cpositivo}(ii)),
which yields $\mM_{\ep_*}^m=\mM_{\ep_*}$. The nonexistence of a third minimal set, the
inequalities $\mm_{\ep_*}\le\bar\muk_{\ep_*}\le \muk_{\ep_*}$, and Proposition~\ref{prop:defM}
ensure that $\mM_{\ep_*}$ is also associated to $\bar\muk_{\ep_*}$ and to $\muk_{\ep_*}$
by \eqref{eq:precollision}. Of course, $x\le\muk_{\ep_*}(\w)$ for any $(\w,x)\in\mM_{\ep_*}$.
In particular, $\mm_{\ep_*}(\w)=\bar\muk_{\ep_*}(\w)=\muk_{\ep_*}(\w)$ for all $\w$ in the residual
subset of $\W$ formed by their common continuity points, and hence $\mM_{\ep_*}$ is contained in
$\bigcup_{\w\in\W} \big(\{\w\}\times [\mm_{\ep_*}(\w),\muk_{\ep_*}(\w)]\big)$,
which is a compact $\tau_{\ep_*}$-invariant pinched subset of $\W\times(0,c_+)$.
\smallskip\par
(ii) Since $c_-\le c_+\le -a(\w)/b(\w)$ for all $\w\in\W$ with $b(\w)\ne 0$, we have
$c_-\le s_0$, where $s_0$ is defined in Theorem \ref{th:cbpositivo}(i). This  result proves (ii).

\smallskip\par
(iii) We fix $\ep<0$, assume $c_+<3\,c_-$, and proceed in a similar way to the proof
of Theorem \ref{th:cpositivo}(v): we fix $r\in(c_+/3,c_-)$; we define
$q_\ep(\w,x)$ as the $C^{0,2}(\WR,\R)$ function which coincide with $p_\ep(\w,x)$ for
$x\ge r$ and is given by a second degree polynomial for $x\le r$; we check that
$(\partial^2/\partial x^2)\,q_\ep(\w,x)=-6\,r+2\,c(\w)\le -6\,r+2\,c_+<0$ for all $x\le r$ and
$(\partial^2/\partial x^2)\,q_\ep(\w,x)=-6\,x+2\,c(\w)<-6\,r+2\,c_+<0$ for all $x>r$;
we deduce that $q_\ep$ satisfies all the conditions {\bf c1}-{\bf c4} of \cite[Section 3]{duno4};
we use $q_\ep(\w,c_-)=p_\ep(\w,c_-)>0$ and $q_\ep(\w,\pm s)<0$ for all $\w\in\W$ if $s>c_-$ is large
enough to deduce from Proposition \ref{prop:encaje} and \cite[Theorem 3.3]{duno4}
the existence of exactly one minimal set $\mM_\ep$ strictly above $\W\times\{c_-\}$
for the skew-product flow $\tilde\tau_\ep$ defined by $x'=q_\ep(\wt,x)$, which is an attractive
hyperbolic copy of the base; and we conclude from $\mA_\ep\subset\W\times(c_-,\infty)$ (proved in (ii))
that $\mM_\ep$ is also the unique $\tau_\ep$-minimal set, and from the coincidence of the
Lyapunov exponents of the minimal sets for both flows which are above $\W\times\{c_-\}$
and from Theorem \ref{th:caract-hiper} that $\mA_\ep=\mM_\ep$ is an attractive hyperbolic
copy of the base, i.e., $\mA_\ep=\{\muk_\ep\}$.
The continuity of $(-\infty,\ep_*)\to C(\W,\R),\,\ep\mapsto\muk_\ep$ follows, again, from the
robustness of the existence of hyperbolic copies of the base
and their continuous variation in the uniform topology (see, e.g., \cite[Theorem 2.3]{duno5}).
\par
\smallskip\par
(iv) Let us take $\ep<0$. If $r\ge s_+$, then $p_\ep(\w,r)=r^2(c(\w)-r)+\ep\,(b(\w)\,r+a(\w))<0$,
since $b(\w)\,r+a(\w)\ge b(\w)\,s_++a(\w)\ge 0$ and $c(\w)-r\le c(\w)-s_+\le c_++a(\w)/b(\w)<0$.
According Remark \ref{rm:comparison}.2, $\muk_\ep<s_+$, as asserted.
Our next goal is to check that, if $-\ep$ is large enough, then
all the Lyapunov exponents of $\mA_\ep$ are strictly negative: according to Theorem \ref{th:caract-hiper},
this property proves that $\mA_\ep=\{\muk_\ep\}$ is an attractive hyperbolic copy of the base.
As explained in Section \ref{subsec:Lyap}, it suffices to check that
\[
 \int_\W(p_\ep)_x(\w,\mb_\ep(\w))\,dm=\int_\W(-3(\mb_\ep(\w))^2+2\,\mb_\ep(\w)\,c(\w)+\ep\,b(\w))\,dm<0
\]
for any $m$-measurable equilibrium $\mb_\ep\colon\W\to\R$ with graph contained in $\mA_\ep$. Since
$c_-\le \mb_\ep\le s_+$ for any such $\tau_\ep$-equilibrium $\mb_\ep$, and since $b_->0$,
the inequality holds if $-\ep>0$ is large enough, say $\ep<\bar\ep<0$. On the other hand,
according to Theorem \ref{th:cpositivo}, $\muk_0\le c_+$ (and hence $\muk_0<s_-$)
and if, in addition, $c_+<3\,c_-$, then $\ep\mapsto\muk_\ep\in C(\W,\R)$ is continuous at $\ep=0$.
So, there exists $\underline\ep\le 0$ (perhaps $\underline\ep=-\infty$)
such that, if $\ep\in(\underline\ep,0)$, then $\muk_{\ep}(\w)\le s_-$ for all $\w\in\W$.
It is easy to deduce that, if $\underline\ep<\ep_1<\ep_2\le 0$,
then $\muk_{\ep_2}'(\w)=p_{\ep_2}(\w,\muk_{\ep_2}(\w))<p_{\ep_1}(\w,\muk_{\ep_2}(\w))$ for all $\w\in\W$, and hence
$\muk_{\ep_2}<\muk_{\ep_1}$ (see once more Remark \ref{rm:comparison}.2).
\smallskip\par
(v) It is clear that $\ep_*$ does not belong to $\mJ$, since according to (i) there are exactly
two $\tau_{\ep_*}$-minimal sets. This ensures that $\ep_*\le\ep^*$. In addition, there exist at most two
$\tau_{\ep^*}$-minimal sets, since $\ep^*\notin J$. We will check below that
there are indeed two of them: $\{\ml_{\ep^*}\}$ (which is hyperbolic attractive and below $\W\times\{0\}$,
as proved by Theorem \ref{th:cpositivo}(ii)) and $\mM_{\ep^*}$, which is above $\W\times\{s_-\}$ for
$s_-:=\inf_{\w\in\W}(-a(\w)/b(\w))>c_+$. Since, as seen in (i), $\mM_{\ep_*}$ is below $\W\times\{c_+\}$,
we conclude that $\ep_*<\ep^*$. The proof will also show that $\mm_\ep>s_-$ for all $\ep\in\mJ$,
and this completes the list of assertions in (v).
\par
It is easy to check that $p_\ep(\w,s_-)<0$ for all $\w\in\W$ if $\ep>0$. So, $s_-$ is a global strict
upper solution for $\tau_\ep$. On the other hand, since $\lim_{\ep\to \infty}\muk_\ep(\w)=\infty$
uniformly on $\W$ (see Theorem \ref{th:cbpositivo}(iii)), there exists a minimum $\ep_1\ge\ep^*$
such that $\muk_\ep(\w)\ge s_-$ for all $\w\in\W$ if $\ep>\ep_1$.
These properties, combined with $\ml_\ep<0<s_-$ for all $\ep>0$ (see Theorem \ref{th:cpositivo}(ii))
mean that, for any fixed
$\w\in\W$ and all $\ep\ge\ep_1$, the map $t\mapsto v_\ep(t,\w,s_-)$ is globally bounded, and that
$\W\times\{s_-\}$ is always strictly above the \upomeg-limit set of $(\w,s_-)$ and
strictly below its \upalfa-limit set: see Remark \ref{rm:comparison}.3.
Remark \ref{rm:comp2} ensures that the \upalfa-limit set
cannot contain $\{\muk_\ep\}$. Hence it necessarily contains $\{\mm_\ep\}$, and the \upomeg-limit
set contains the unique $\tau_\ep$-minimal set below $\{\mm_\ep\}$, which is $\{\ml_\ep\}$.
In particular, $\mm_\ep>s_-$ for all $\ep>\ep_1$.
Let us check that $\ep_1=\ep^*$, assuming for contradiction that $\ep_1>\ep^*$. Using the continuity of
$\ep\mapsto\muk_\ep$ on $\mJ$ ensured by Theorem \ref{th:cbpositivo}(iii), we deduce that $\muk_{\ep_1}(\w)\ge s_-$
for all $\w\in\W$, and that there exists $\w_0\in\W$ such that $\muk_{\ep_1}(\w_0)=s_-$. But, as just seen,
this yields $\muk_{\ep_1}>\mm_{\ep_1}>s_-$, which provides the contradiction. Altogether, we
have $\mm_\ep>s_->\ml_\ep$ for all $\ep\in\mJ$.
\par
Note that we have proved that $\W\times\{s_-\}\subset\mA_\ep$
for all $\ep>\ep^*$. Since there exists $\rho>0$ such that $s_-<\mm_\ep<\muk_\ep\le\rho$ for all $\ep\in(\ep^*,\ep^*+1]$,
we conclude that $v_{\ep^*}(t,\w,s_-)=\lim_{\ep\to(\ep^*)^+}v_\ep(t,\w,s_-)\le\rho$ for all $t<0$ and $\w\in\W$.
So, $\ml_{\ep^*}<s_-\le \muk_{\ep^*}$, and hence, as seen above,
there exists the \upalfa-limit set of any point $(\w,s_-)$ for $\tau_{\ep^*}$ and it is strictly above
$\W\times\{s_-\}$. Hence, there exists a $\tau_{\ep^*}$-minimal set contained in $\W\times(s_-,\infty)$,
as asserted.
\smallskip\par
(vi) Assume that $-a(\w)/b(\w)=s$, constant.
Let us first analyze the situation in $(\ep^*,\infty)$. Recall that $\muk_\ep>\mm_\ep>s$
if $\ep>\ep^*$: see (v). An argument similar to that of the first (resp.~second) paragraph
of the proof of (i) proves that $\ep\mapsto\muk_\ep$ (resp.~$\ep\mapsto-\mm_\ep$) is strictly increasing
on $(\ep^*,\infty)$. In particular, there exist the pointwise limits $\bar\muk_{\ep^*}:=\lim_{\ep\to(\ep^*)^+}\muk_\ep$
and $\mm_{\ep^*}:=\lim_{\ep\to(\ep^*)^+}\mm_\ep$, whose additional properties are checked as those of
$\bar\muk_{\ep_*}$ and $\mm_{\ep_*}$ in (i). Let us prove that $\lim_{\ep\to\infty}\mm_\ep(\w)=s$ uniformly
on $\W$. We take $\delta>0$ such that, for an $\ep_2>\ep^*$, $\ml_\ep<0<s+\delta<\muk_\ep$ if $\ep>\ep_2$,
and work for these values of $\ep$. Note that $p_\ep(\w,s+\delta)=(s+\delta)^2(c(\w)-s-\delta)+\ep\,b(\w)\,\delta$,
so that $s+\delta$ is a global strict lower solution if $\ep>0$ is large enough.
According to Remark \ref{rm:comparison}.3,
$\W\times\{s+\delta\}$ is strictly above a $\tau_\ep$-minimal set
contained in the \upalfa-limit set of a point $(\w,s+\delta)$ (that exists since
this point belongs to $\mA_\ep$), which is necessarily $\{\mm_\ep\}$,
as Remark \ref{rm:comp2} guarantees.
That is, $s\le\mm_\ep\le s+\delta$ if $\ep$ is large enough, and this proves the assertion.
\par
Let us now analyze the situation for $\ep\in(-\infty,\ep_*)$. First of all, we check that
$\lim_{\ep\to-\infty}\muk_\ep(\w)=s$ uniformly on $\W$. For $\delta>0$, $p_\ep(\w,s-\delta)=
(s-\delta)^2(c(\w)-s+\delta)-\ep\,b(\w)\,\delta$, so $s-\delta$ is a global strict lower solution if $-\ep>0$ is
large enough. According to Remark \ref{rm:comparison}.3, $\W\times\{s-\delta\}$ is strictly below a $\tau_\ep$-minimal set
contained in the \upomeg-limit set of a point $(\w,s-\delta)$; and, according to (iv), this $\tau_\ep$-minimal set is necessarily
$\{\muk_\ep\}$ if $\ep<\bar\ep$.
That is, $s-\delta<\muk_\ep<s$ if $-\ep>0$ is large enough (see (iv)), which proves the assertion.
To check that $\ep\mapsto\muk_\ep$ is strictly decreasing on $(-\infty,\ep_*)$, we repeat the argument of (iv),
which is possible since, by (iv) and (i), $\muk_\ep<s$ for all these values of $\ep$. (Recall that, in
addition, $\mA_\ep=\{\muk_\ep\}$ for all $\ep\le 0$ if $c_+<3\,c_-$: see (iii).)
\par
Let us see what happens for $\ep\in(\ep_*,\ep^*)$. We fix $\bar\ep\in(\ep_*,\ep^*)$.
As a first step, we will check that $\mA_{\bar\ep}\subset\W\times(-\infty,s)$; i.e.,
that $s>\muk_{\bar\ep}(\w)$ for all $\w\in\W$. Since $\ml_{\bar\ep}<0$ (see Theorem \ref{th:cpositivo}(ii)),
it suffices to assume that the $\tau_{\bar\ep}$-orbit of a point $(\w_0,s)$ is bounded (i.e.,
that $s\le\muk_{\bar\ep}(\w)$) and reach a contradiction.
Since $p_{\bar\ep}(\w,s)=s^2\,(c(\w)-s)<0$, $s$ is a constant strict upper solution
for $\tau_{\bar\ep}$. Then, as seen in the proof of (v), $\W\times\{s\}$ is strictly above
the \upomeg-limit set of $(\w_0,s)$, which ensures that $s<\muk_{\bar\ep}$. So,
we can repeat once again the argument of (iv) to check that $\ep\mapsto\muk_\ep$ is strictly increasing on $(\bar\ep,\infty)$.
In particular, the $\tau_\ep$-orbit of $(\w_0,s)$ is bounded for all $\ep\ge\bar\ep$.
Let $\bar\w$ be a continuity point of the semicontinuous map $\muk_{\ep^*}$, so that
$(\mM_{\ep^*})(\bar\w)=\{\muk_{\ep^*}(\bar\w)\}$ (see Proposition \ref{prop:defM}).
The \upalfa-limit set of $(\w_0,s)$ for $\tau_{\ep^*}$ contains a minimal set which cannot be
hyperbolic attractive (see again Remark \ref{rm:comp2}), so it is $\mM_{\ep^*}$.
That is, there exists $(t_n)\downarrow-\infty$ such that
$(\bar\w,\muk_{\ep^*}(\bar\w))=\lim_{n\to\infty}(\w_0{\cdot}t_n,v_{\ep^*}(t_n,\w_0,s))$.
We can assume without restriction the existence of $\bar x:=\lim_{n\to\infty} v_{\bar\ep}(t_n,\w_0,s)$, and
observe that $\bar x\le\muk_{\bar\ep}(\bar\w)$, since $v_{\bar\ep}(t,\bar\w,\bar x)$ is bounded
($v_{\bar\ep}(t,\bar\w,\bar x)=\lim_{n\to\infty} v_{\bar\ep}(t_n+t,\w_0,s)$). Then, since
$v_{\bar\ep}(t,\w_0,s)>v_{\ep^*}(t,\w_0,s)$ for all $t<0$ (as we deduce from $p_{\bar\ep}(\w,r)<p_{\ep^*}(\w,r)$ for all
$\w\in\W$ if $r\ge s$ and from $v_\ep(t,\w_0,s)>s$ for all $t<0$ and all $\ep>0$),
we can conclude that $\bar x\ge\muk_{\ep^*}(\bar\w)>\muk_{\bar\ep}(\bar\w)$,
which provides the sought-for contradiction.
\par
So, we have $\mA_{\bar\ep}\subset\W\times(-\infty,s)$. Now we assume for contradiction the
existence of a $\tau_{\bar\ep}$-minimal set $\mN_{\bar\ep}>\{\ml_{\bar\ep}\}$.
We choose a point $\bar\w\in\W$ at which the sections of $\mM_{\ep_*}$ and $\mN_{\bar\ep}$ are singletons:
$(\mM_{\ep_*})_{\bar\w_1}=\{\muk_{\ep_*}(\bar\w)\}=\{\mm_{\ep_*}(\bar\w)\}$ and
$(\mN_{\bar\ep})_{\bar\w}=\{\bar x\}$, so that $\bar x\in (\ml_\ep(\bar\w),s)$.
The information regarding monotonicity, continuity, limiting
behavior as $\ep\to-\infty$, and the shape of $\mM_{\ep_*}$ provided so far, allows us to choose
an $\ep_0<\bar\ep$ and a unique (bounded) $\tau_{\ep}$-equilibrium, say $\mb_{\ep_0}$,
such that $\mb_{\ep_0}(\bar \w)=\bar x$: if $\bar x\in(\ml_\ep(\bar\w),0]$,
then $\ep_0\in[0,\bar\ep)$ and $\mb_{\ep_0}=\ml_{\ep_0}$; if $\bar x\in(0,\mm_{\ep_*}(\bar\w))=
(0,\muk_{\ep_*}(\bar\w))$, then $\ep\in(0,\ep_*)$ and $\mb_{\ep_0}=\mm_{\ep_0}$; and if
$\bar x\in[\muk_{\ep_*}(\bar\w),s)$, then $\ep_0\in(-\infty,\ep_*]$ and
$\mb_{\ep_0}=\muk_{\ep_0}$. In any case, $\mb_{\ep_0}$ is a global strict upper solution for $\tau_{\bar\ep}$,
since $\mb_{\ep_0}'(\w)=p_{\ep_0}(\w,\mb_{\ep_0}(\w))>p_{\bar\ep}(\w,\mb_{\ep_0}(\w))$ due to the inequality
$\mb_\ep<s$. As explained in Remark \ref{rm:comparison}.3,
this ensures that the \upomeg-limit of the point
$(\bar\w,\mb_{\ep_0}(\bar\w))=(\bar\w,\bar x)$
(which is, of course, $\mN_{\bar\ep}$) is strictly below the graph of $\mb_{\ep_0}$.
This fact precludes
$(\bar\w,\bar x)\in\mN_{\bar\ep}$ and provides the sought-for contradiction.
\par
The assertions of the last sentence of (vi) follow from the previous description and Theorems \ref{th:cbpositivo}(i) and
\ref{th:cpositivo}(v).
\end{proof}
Figure \ref{fig:caso1} provides a depiction of the ``three saddle-node
bifurcation diagram" of \eqref{eq:EDOpol} under the most restrictive conditions of
Theorem \ref{th:caso1}. It is interesting to remark that the dynamics of the nonhyperbolic
minimal sets $\mM_{\ep_*}$ and $\mM_{\ep^*}$ at the bifurcation points $\ep_*$ and $\ep^*$
can be extremely complicated, even with the occurrence of SNAs described, for instance, in
\cite{jnot}, \cite{nuob8}, \cite{jage, jage2} (based on the
classical examples of \cite{mill,mill2}, \cite{vino} and \cite{john14}). A more detailed description
of these dynamical possibilities can be easily adapted to this case from that made in
\cite[Proposition 5.11]{duno1}. In particular, as there explained, the measure $m$ of the residual
subsets $\mR^{\ep_*}$ and $\mR^{\ep^*}$ of $\W$ at whose points the upper and lower equilibria of
$\mM_{\ep_*}$ and $\mM_{\ep^*}$ respectively collide can be $0$ or $1$.
\begin{figure}[h]
\centering
\includegraphics[width=0.9\textwidth]{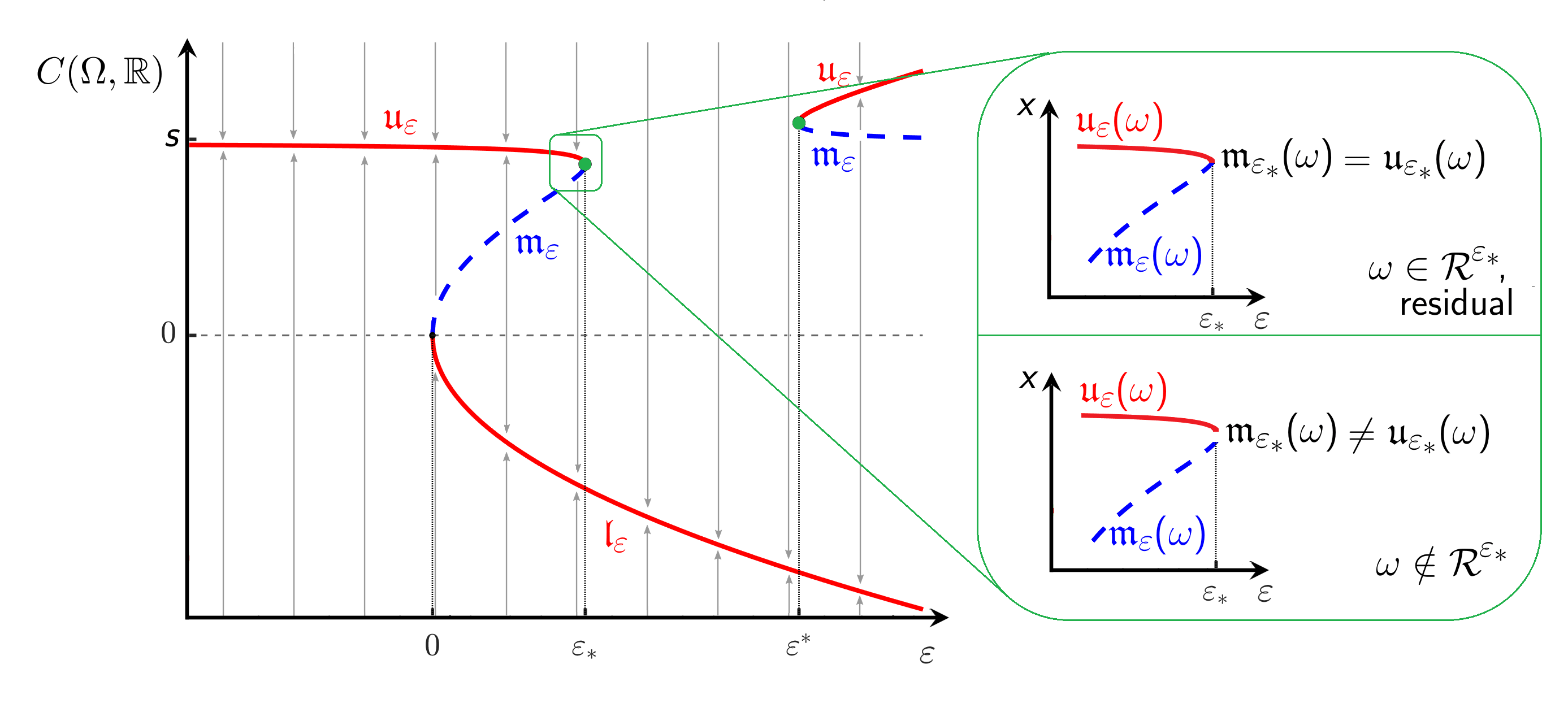}
\caption{The figure in the left depicts the bifurcation diagram of the $\ep$-parametric family \eqref{eq:EDOpol}
when $c>0$, $c_+<3c_-$, $b>0$, $a=-s\,b$ for a constant $s>0$, and $b\,c_++a<0$,
which is described in Theorem  \ref{th:caso1}
in combination with Theorems \ref{th:cpositivo} and \ref{th:cbpositivo}.
The three strictly monotone solid red curves represent the families of attractive hyperbolic copies of the base that
determine the upper and lower equilibria of the global attractor. The two strictly monotone dashed blue curves
represent the families of repulsive hyperbolic copies of the base. The light grey arrows partly
depict the dynamics of the rest of the orbits. There are three bifurcation
points, all of them of saddle-node type: $\ep=0$, where $\mm_\ep$ and $\ml_\ep$ globally collide on $0$;
and $\ep_*$ and $\ep^*$, where $\mm_\ep$ and $\muk_\ep$ partly collide, giving rise to semicontinuous but
perhaps noncontinuous maps. This fact is depicted by two large green points, and explained for $\ep_*$ in the
zoom at the right. (Figures \ref{fig:casos23} and \ref{fig:caso4} will also use ``large green points"
to depict similar situations, as well as red and blue curves, and grey arrows.)}
\label{fig:caso1}
\end{figure}
%
\begin{teor}\label{th:caso2}
Assume that $a(\w)<0$, $b(\w)>0$, $c(\w)>0$ and $b(\w)\,c_-+a(\w)>0$ for all $\w\in\W$, and
call $s_+:=\sup_{\w\in\W}(-a(\w)/b(\w))$ and $s_-:=\inf_{\w\in\W}(-a(\w)/b(\w))$.
Then, in addition to the information provided by Theorems {\rm\ref{th:cpositivo}} and {\rm\ref{th:cbpositivo}},
\begin{itemize}[leftmargin=20pt]
\item[\rm(i)] for all $\ep>0$, there are three hyperbolic copies of the base, with $\ml_\ep<0<\mm_\ep<s_+<\muk_\ep$. That is,
    $\mI=(0,\infty)$, where $\mI$ is defined in Theorem {\rm\ref{th:cpositivo}(iii)}.
    In addition,
    the maps $(0,\infty)\to C(\W,\R),\,\ep\mapsto-\ml_\ep,\muk_\ep$ are strictly increasing; and there exists
    $\ep^*$ with $0<\ep^*\le\infty$ such that the map $(0,\ep^*)\to C(\W,\R),\,\ep\mapsto\mm_\ep$ is
    strictly increasing. In particular, there are no strictly positive bifurcation values.
\item[\rm(ii)] For all $\ep<0$, $\mA_\ep\subset\W\times(s_-,c_+)$; there exists $\ep_*$ with $-\infty\le\ep_*<0$
    such that the map $(\ep_*,\infty)\to C(\W,(s_-,c_+)),\,\ep\mapsto\muk_\ep$ is strictly increasing;
    and there exists $\ep_0\le 0$ such that $\mA_\ep$ is an attractive hyperbolic copy of the base
    for $\ep\in(-\infty,\ep_0)$.
\item[\rm(iii)] If, in addition, $c_+<3\,s_-$, then $\ep_0=0$. Hence,
    $\{\muk_\ep\}$ is an attractive hyperbolic copy of the base for all $\ep\in\R$ and
    the map $\R\to C(\W,\R),\,\ep\mapsto\muk_\ep$ is continuous. Consequently, in this case,
    $0$ is the unique bifurcation value, of local saddle-node type.
\item[\rm(iv)] If, in addition, $a/b=-s\in\R$, then $\ep^*=\infty$ and $\ep_*=-\infty$;
    the map $(-\infty,0)\to C(\W,(s,c_+)),\,\ep\mapsto\ml_\ep$ is strictly increasing; and
    $\lim_{\ep\to-\infty}\muk_\ep(\w)=\lim_{\ep\to\infty}\mm_\ep(\w)=s$ uniformly on $\W$.
    There are three possibilities for $\ep\in(-\infty,0)$:
    \begin{itemize}[leftmargin=10pt]
    \item[-] $\mA_\ep=\{\muk_\ep\}$ is an attractive hyperbolic copy of the base for all $\ep<0$, in which case
    the map $\R\to C(\W,\R),\,\ep\mapsto\muk_\ep$ is continuous, and $0$ is the
    unique bifurcation value, of local saddle-node type. (This happens if $c_+<3\,s$.)
    \item[-] There exist $\underline\ep<\bar\ep<0$ such that: there are
    three $\tau_\ep$ hyperbolic copies of the base for any $\ep\in(\underline\ep,\bar\ep)$;
    $\mA_\ep$ is an attractive hyperbolic copy of the base for $\ep\in(0,\infty)-[\,\underline\ep,\,,\,\bar\ep\,]$;
    there are two $\tau_{\underline\ep}$-minimal sets, $\{\ml_{\underline\ep}\}$
    (which is hyperbolic attractive) and a nonhyperbolic one given by the collision of $\{\muk_\ep\}$ and $\{\mm_\ep\}$
    as $\ep\to(\underline\ep)^+$; there are two $\tau_{\bar\ep}$-minimal sets, $\{\muk_{\bar\ep}\}$
    (which is hyperbolic attractive) and a nonhyperbolic one given by the collision of $\{\ml_\ep\}$ and $\{\mm_\ep\}$
    as $\ep\to(\bar\ep)^-$; the maps $(-\infty,\bar\ep)\to C(\W,\R),\,\ep\to\ml_\ep$ and $(\underline\ep,0)\to C(\W,\R),\,\ep\to\muk_\ep$
    are continuous and strictly increasing; and the map $(\underline\ep,\bar\ep)\to C(\W,\R),\,\ep\to\mm_\ep$
    is continuous and strictly decreasing. So, there are exactly three bifurcation values,
    $\underline\ep$, $\bar\ep$ and $0$, all of them of local saddle-node type.
    \item[-] There is a unique negative value $\ep_1<0$ such that $\{\muk_{\ep_1}\}$ is not a hyperbolic copy of the base,
    in which case $\mA_{\ep_1}$ is a pinched set containing a unique $\tau_{\ep_1}$-minimal set.
    \end{itemize}
\end{itemize}
\end{teor}
\begin{proof}
(i) If $\ep\ge0$, then $p_\ep(\w,s_+)=(s_+)^2(-s_++c(\w))+\ep\,(b(\w)\,s_++a(\w))>0$, since
$b(\w)\,s_++a(\w)\ge 0$ for all $\w\in\W$ and $c_->-a(\w)/b(\w)$ for all $\w$ (and hence $c(\w)\ge c_->s_+$).
In addition, $p_\ep(\w,0)=\ep\,a(\w)<0$ if $\ep>0$. We fix $\ep>0$, look for $r_1<0<s_+<r_2$ with $p_\ep(\w,r_1)>0$
and $p_\ep(\w,r_2)<0$ for all $\w\in\R$, and deduce the first assertion in (i) from
Proposition \ref{prop:encaje} and Theorem \ref{th:nummin}. If $0<\ep_1<\ep_2$, then $\muk_{\ep_1}'(\w)<
p_{\ep_2}(\w,\muk_{\ep_1}(\w))$ (since $b(\w)\,\muk_{\ep_1}(\w)+a(\w)>b(\w)\,s_++a(\w)\ge 0$),
and hence, as explained in Remark \ref{rm:comparison}.2, $\muk_{\ep_1}<\muk_{\ep_2}$: $\ep\mapsto\muk_\ep$ is
strictly increasing on $(0,\infty)$.
To prove the monotonicity of $\ep\mapsto\mm_\ep$, we use an argument analogous to that of the
last paragraph in the proof of Theorem \ref{th:cpositivo}(iii), working on the interval $(0,\ep^*)$
(perhaps $(0,\infty)$) on which $b(\w)\,\mm_\ep(\w)+a(\w)<0$: this happens if $\mm_\ep<s_-$, and hence
at least for small values of $\ep>0$, since $\lim_{\ep\to 0^+}\mm_\ep(\w)=0$ uniformly on $\W$.
The lack of strictly positive bifurcation values of $\ep$ is a trivial consequence of the previous properties.
\smallskip\par
(ii) Let us fix $\ep<0$. By hypotheses, $s_-=s_0$, with $s_0$ defined in Theorem \ref{th:cbpositivo}(i), and
this results ensures that $\mA_\ep\subset\W\times(s_-,\infty)$. In addition, if $r\ge c_+$,
then $p_\ep(\w,r)<0$ for all $\w\in\W$, since $c(\w)-r\le 0$ and $b(\w)\,r+a(\w)\ge b(\w)\,c_-+a(\w)>0$.
According to Remark \ref{rm:comparison}.1, $\mA_\ep\subset\W\times(s_-,c_+)$.
\par
Since $\{\muk_0\}$ is a hyperbolic copy of the base and $\muk_0\ge c_->s_+$ (see Theorem \ref{th:cpositivo}(i)),
the persistence ensured by \cite[Theorem 2.3]{duno5} guarantees the existence of a hyperbolic
copy of the base strictly greater than $s_+$ for $\ep<0$ close enough to $0$, and hence the definition
of $\muk_\ep$ shows that $\muk_\ep>s_+$ for these values of $\ep$.
Let $(\ep_*,0)\subseteq(-\infty,0)$ be the interval of persistence of
this property (on which we cannot guarantee the continuity of $\muk_\ep$).
If $\ep_*<\ep_1<\ep_2<0$, then $\muk_{\ep_1}'(\w)<p_{\ep_2}(\w,\muk_{\ep_1}(\w))$
(since $b(\w)\,\muk_{\ep_1}(\w)+a(\w)>b(\w)\,s_++a(\w)\ge 0$), and
hence Remark \ref{rm:comparison}.2 ensures that $\muk_{\ep_1}<\muk_{\ep_2}$, as asserted.
\par
Let us check that $\mA_\ep$ is an attractive hyperbolic copy of the base if $-\ep$ is large enough.
We deduce from $b_->0$ that $\int_\W(p_\ep)_x(\w,\mb_\ep(\w))\,dm<0$ for all $m$-measurable
equilibrium $\mb_\ep\colon\W\to(s_-,c^+)$ (i.e., with graph contained in $\mA_\ep$) if, let's say,
$\ep<\ep_0<0$. As explained in the proof of Theorem \ref{th:caso1}(iv), this ensures
that all the Lyapunov exponents of the global attractor are strictly negative, and hence
the assertion follows from Theorem \ref{th:caract-hiper}.
\smallskip\par
(iii) Let us assume that $c_+<3\,s_-$.
To check the first assertion in (iii), we use again an argument similar
to that used in the proofs of Theorems \ref{th:cpositivo}(v) and \ref{th:caso1}(iii).
We fix $\ep<0$; we define $q_\ep(\w,x)$ as the $C^{0,2}(\WR,\R)$ function
which coincides with $p_\ep(\w,x)$ for $x\ge s_-$ and is given by a second degree
polynomial for $x\le s_-$; we check that $(\partial^2/\partial x^2)\,q_\ep(\w,x)<0$ for
$x\ge s_-$ and for $x<s_-$; and we deduce from this and its shape that
$q_\ep$ satisfies all the conditions {\bf c1}-{\bf c4} of \cite[Section 3]{duno4}.
In addition, $s_-$ is a strict global lower solution for our $\ep<0$, since $q_\ep(\w,s_-)=p_\ep(\w,s_-)>0$
(which follows from $s_-<c_+$). Hence, there exists exactly a minimal
set $\mM^u_\ep$ for the skew-product flow $\tilde\tau_\ep$ defined by $x'=q_\ep(\wt,x)$,
and it is hyperbolic attractive and strictly above $\W\times\{s_-\}$. Since any
$\tau_\ep$-minimal set is strictly above $\W\times\{s_-\}$, $\mM^u_\ep$ is the unique one,
and since its Lyapunov exponents are the same for $p_\ep$ as for $q_\ep$ (i.e., negative), then
Theorem \ref{th:caract-hiper} ensures that $\mA_\ep=\mM_\ep$ is an attractive hyperbolic copy of the base,
$\mA_\ep=\{\muk_\ep\}$. The usual persistence argument shows the continuity of $\ep\mapsto\muk_\ep$
on $\R$, and the last assertion in (iii) is a consequence of the previous analysis.
\smallskip\par
(iv) Let us assume that $a/b=-s\in\R$. By reviewing the proofs of (i) and (ii), we see that $\ep^*=\infty$
(since $s_-=s_+$ and hence $\mm_\ep<s_-$ for all $\ep>0$) and $\ep_*=-\infty$ (since $s_+=s_-$ and hence
$\muk_\ep>s_+$ for all $\ep<0$).
We fix $\ep_1<\ep_2<0$ and deduce from $\ml_{\ep_2}>s_-=s_+$ (see (ii)) that
$\ml_{\ep_2}'(\w)>p_{\ep_1}(\w,\ml_{\ep_2}(\w))$ for all $\w\in\W$, so $\ml_{\ep_2}(\w)>\ml_{\ep_1}(\w)$
(see Remark \ref{rm:comparison}.2). Let us check that
$\lim_{\ep\to\infty}\mm_\ep(\w)=s$ uniformly on $\W$. We take $\delta>0$. Since $p_\ep(\w,s-\delta)=
(s-\delta)^2(c(\w)-s+\delta)-\ep\,b(\w)\,\delta$, we have that
$s-\delta$ is a global strict upper solution if $\ep>0$ is
large enough, and $(\w,s-\delta)\in\mA_\ep$ for all $\w\in\W$ and a large enough $\ep$ (see (i)),
Remark \ref{rm:comparison}.3 shows that, for these values of $\ep$,
$\W\times\{s-\delta\}$ is strictly below a $\tau_\ep$-minimal set
contained in the \upalfa-limit set of a point $(\w,s-\delta)$, which according to Remark \ref{rm:comp2}
is necessarily $\{\mm_\ep\}$.
That is, $\mm_\ep>s-\delta$ if $\ep$ is large enough, which combined with $\mm_\ep<s$ proves the assertion. A
similar argument, working with $s+\delta$, shows that $\lim_{\ep\to-\infty}\muk_\ep(\w)=s$ uniformly on $\W$.
\par
The property $\mA_\ep\subset\W\times(s,c_+)$ for any $\ep<0$ (see (iii))
allows us to check that, if $\ep_1<\ep_2<0$, any $\tau_{\ep_1}$-equilibrium is a global
strict lower solution for $\tau_{\ep_2}$, and any $\tau_{\ep_2}$-equilibrium is a
global strict upper solution for $\tau_{\ep_1}$. These properties are required several times
in the steps leading to a detailed proof of the remaining assertions, which we only sketch.

Let us assume the existence of an $\ep_1<0$ such that $\mA_{\ep_1}$ is not a hyperbolic copy of the base
(which, according to (iii), ensures that $\ep_0\le\ep_1$ and is not possible if $c_+<3\,s$). To start with,
we also assume that there are three $\tau_{\ep_1}$-minimal sets. We call $\mI$ the maximal interval containing
$\ep_1$ at which this property holds. We know that $\mI\subset[\ep_0,0]$ (see point (ii) and Theorem
\ref{th:cpositivo}(i)), and that it is open. This property and those mentioned in the previous
paragraph allow us to repeat the arguments leading to the proof of \cite[Theorem 5.10]{duno1}
to conclude the existence of $\underline\ep$ and $\bar\ep$ with $\ep_0\le\underline\ep<\ep_1<\bar\ep\le 0$
satisfying the stated properties. To check that $\bar\ep<0$, it is enough to observe that the lower
$\tau_{\bar\ep}$-minimal set is strictly above $\{0\}$, which is the lower $\tau_0$-minimal set
(see again Theorem \ref{th:cpositivo}(i)).
\par
Now we assume that there are exactly two $\tau_{\ep_1}$-minimal sets. According to \cite[Theorem 5.13(iii)]{duno1},
one of them is hyperbolic attractive. This allows us to repeat the arguments of \cite[Theorem 5.12]{duno1}
to conclude the existence of three $\tau_{\ep_2}$-minimal sets for an $\ep_2<0$ close to $\ep_1$,
and hence we are in the same situation of the previous paragraph (being in this case $\ep_1$ one of the
two negative bifurcation values).
\par
The remaining case is that $\mA_{\ep_1}$, which is not a hyperbolic copy of the base, contains just
one $\tau_{\ep_1}$-minimal set, which is necessarily nonhyperbolic: see Theorem \ref{th:caract-hiper}(iii).
The previous analysis
shows that there exists just one $\tau_\ep$-minimal set for any $\ep<0$, and hence we can reason as in
\cite[Theorems 5.14 and 5.15]{duno1} to conclude that the situation is the last one described in the statement
of the theorem.
\end{proof}
It is easy to find autonomous examples of the three cases described in the previous point (iv),
what makes sense of this case study. Figure \ref{fig:casos23} depicts two of these three
bifurcation diagrams of \eqref{eq:EDOpol}, appearing under the most restrictive conditions of
Theorem \ref{th:caso2}.
\begin{figure}[h]
\centering
\includegraphics[width=\textwidth]{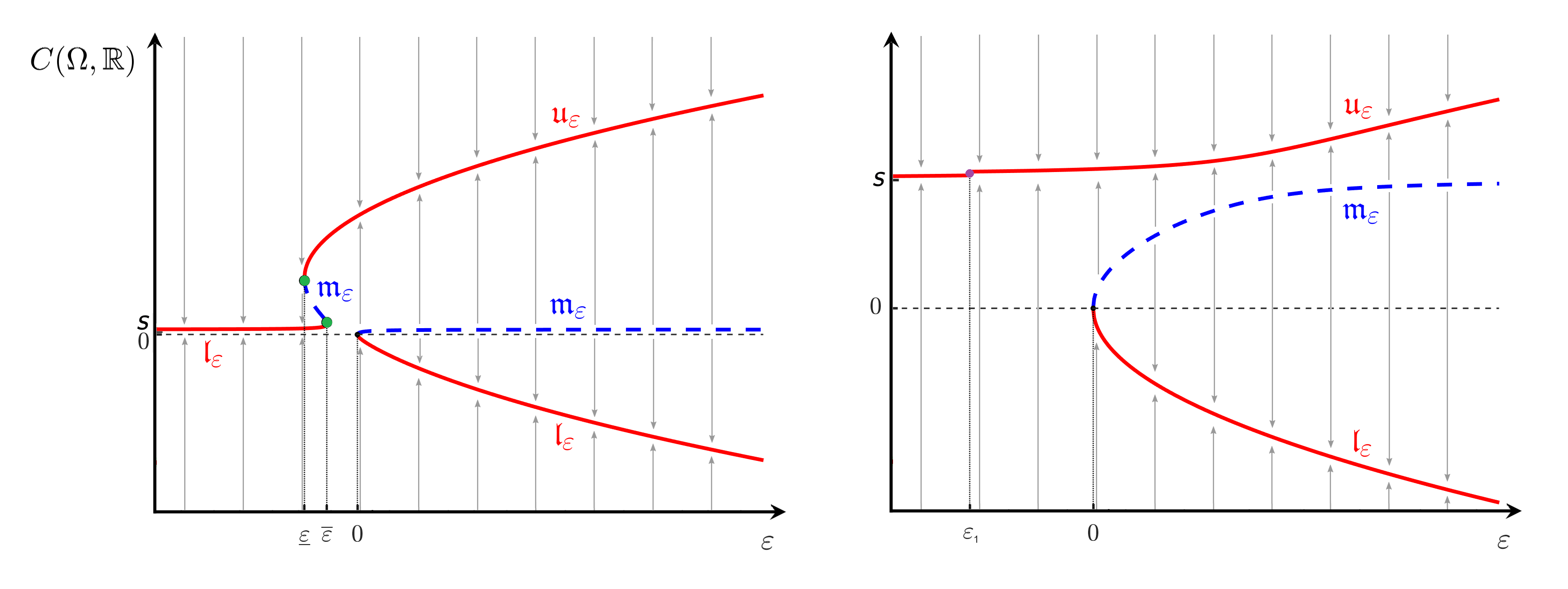}
\caption{The two panels represent two of the three possibilities for the bifurcation diagram of
the $\ep$-parametric family \eqref{eq:EDOpol} when
when $c>0$, $b>0$, $a=-s\,b$ for $s\in(0,\infty)$, and $s<c$.
See Theorem \ref{th:caso2} (in combination with Theorems \ref{th:cpositivo} and \ref{th:cbpositivo})
for the results, and
the caption of Figure \ref{fig:caso1} for the meaning of the different elements.
In the left panel, there are three bifurcation values of the parameter, all of them of local saddle-node type.
In the right panel,
the purple point over $\ep_1$ depicts a pinched attractor containing a unique nonhyperbolic $\tau_{\ep_1}$-minimal
set. In the non depicted bifurcation diagram, which holds at least if, in addition,
$c<3s$, a solid red upper continuous curve would represent the evolution of $\muk_{\ep}$ as $\ep$ varies in $\R$, and hence
$\ep=0$ would be the unique bifurcation value, of local saddle-node type.
}
\label{fig:casos23}
\end{figure}

\begin{figure}[h]
\centering
\includegraphics[width=0.55\textwidth]{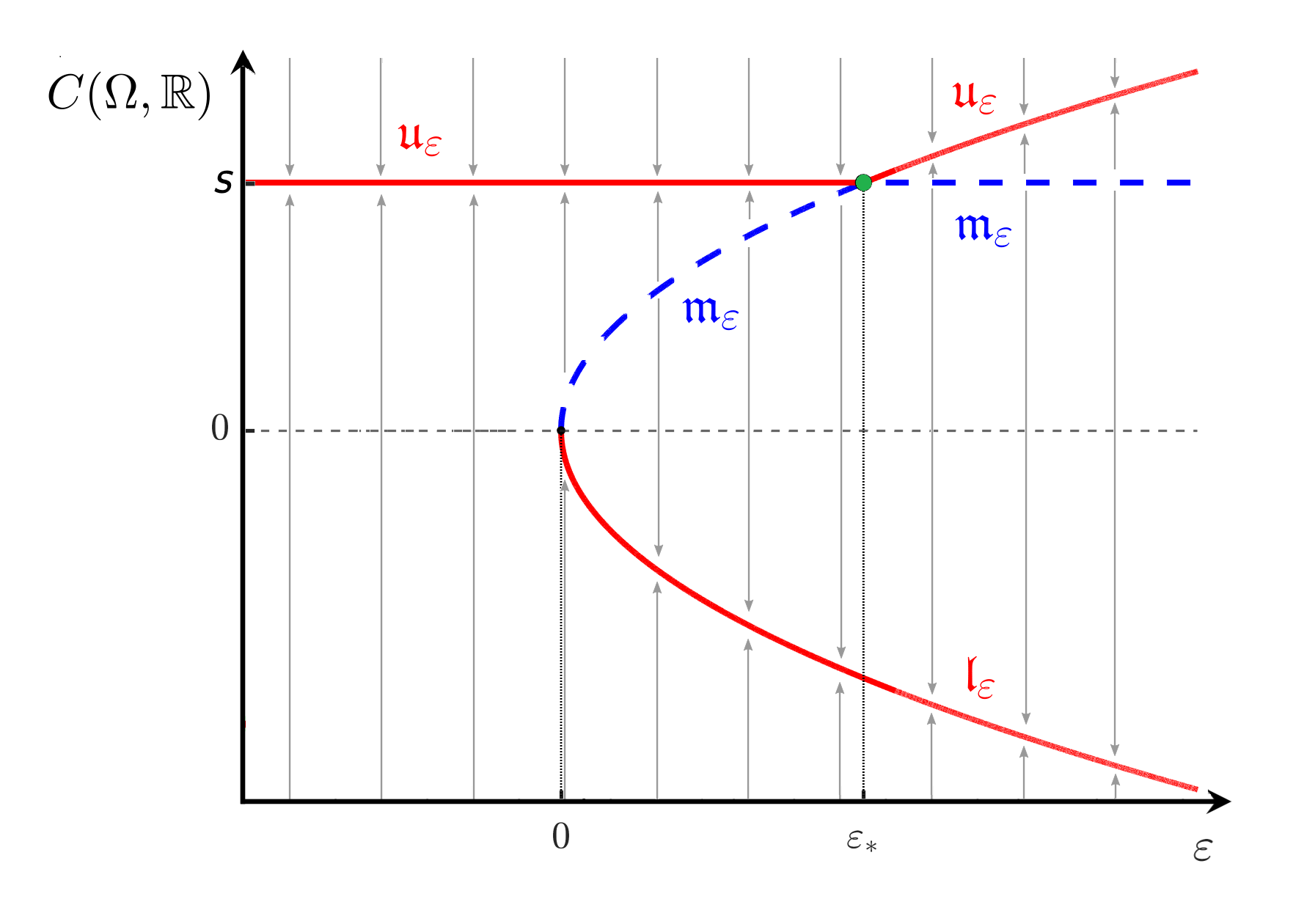}
\caption{The bifurcation diagram of the $\ep$-parametric family \eqref{eq:EDOpol}
when $c=s$ for $s\in(0,\infty)$, $b>0$ and $a=-s\,b$. In this case, $0$ is a local
saddle-node bifurcation point, and $\ep^*:=s^2/\int_\W b(\w)\,dm$ is a transcritical bifurcation point.
This is proved in Theorem \ref{th:caso3}, combined with Theorems \ref{th:cpositivo} and \ref{th:cbpositivo}.
The meaning of the different elements is explained in Figure~\ref{fig:caso1}.}
\label{fig:caso4}
\end{figure}

The hypotheses of the next theorem are considerably more restrictive than those of the previous ones.
We include this result by completeness, and point out that it completes the description of
all the possibilities
for the bifurcation diagrams in the autonomous case with $c>0$, $b>0$ and $a<0$.
Recall that $\ep_0$ is a (nonautonomous) local transcritical bifurcation point
when two hyperbolic copies of the base exist for close values of $\ep$  and approach each
other as $\ep\to(\ep_0)$, giving rise to a locally unique $\tau_{\ep_0}$-minimal set,
which is nonhyperbolic.
\begin{teor}\label{th:caso3}
Assume that $b(\w)>0$, $c(\w)\equiv s>0$ (constant), and $a(\w)=-s\,b(\w)$ for all $\w\in\W$.
Let $\ep_*:=\sup\mI$, with $\mI$ defined in Theorem {\rm\ref{th:cpositivo}(iii)}. Then,
$\ep_*=s^2/\int_\W b(\w)\,dm$. In addition, $\{s\}$ is a $\tau_\ep$-copy of the base for all $\ep\in\R$,
it is hyperbolic attractive with $s=\muk_\ep$ if and only if $\ep<\ep_*$,
and it is hyperbolic repulsive with $s<\muk_\ep$ if and only if $\ep>\ep_*$. More precisely,
in addition to the information provided by Theorems {\rm\ref{th:cpositivo}} and {\rm\ref{th:cbpositivo}},
\begin{itemize}[leftmargin=20pt]
\item[\rm(i)] if $\ep<0$, then $\mA_\ep=\{s\}$, and it is hyperbolic attractive.
\item[\rm(ii)] There exist exactly two $\tau_0$-minimal sets: $\{\ml_0\}=\{0\}$, which is nonhyperbolic, and
    $\{\muk_0\}=\{s\}$, hyperbolic attractive.
\item[\rm(iii)] If $\ep>0$ and $\ep\ne\ep_*$, there are three hyperbolic copies of the base, given by $\ml_\ep<\mm_\ep<s$ for
    $\ep\in(0,\ep_*)$ and by $\ml_\ep<s<\muk_\ep$ for $\ep>\ep_*$.
\item[\rm(iv)] There exist exactly two $\tau_{\ep_*}$-minimal sets: $\{\ml_\ep\}$, which is hyperbolic attractive, and
    $\{s\}$, nonhyperbolic.
\item[\rm(v)] The maps $[0,\infty)\to C(\W,\R),\,\ep\mapsto-\ml_\ep$,
    $[0,\ep_*)\to C(\W,\R),\,\ep\mapsto\mm_\ep$ (with $\mm_0:=0$), and $(\ep_*,\infty)\to C(\W,\R),\,\ep\mapsto\muk_\ep$
    are continuous and strictly increasing.
    In addition,
    the semicontinuous maps $\muk_{\ep_*}$ and
    $\mm_{\ep_*}:=\lim_{\ep\to(\ep^*)^-}\mm_\ep$ take the value $s$ at their continuity points.
\end{itemize}
Therefore, there exist exactly two bifurcation points: $0$, of local saddle-node type, and $\ep_*$, of local transcritical type.
\end{teor}
\begin{proof}
Note that the hypotheses guarantee that $p_\ep(\w,x)=(x-s)\,(-x^2+\ep\,b(\w))$.
Since $p_\ep(\w,s)=0$ for all $\w\in\W$ and $\ep\in\R$, $\{s\}$ is a $\tau_\ep$-copy of the base; and
the remaining initial assertions follow from the fact that $\int_\W p_x(\w,s)\,dm=-s^2+\ep\int_\W b(\w)\,dm$ is its unique
Lyapunov exponent: see Section \ref{subsec:Lyap}.
\smallskip\par
(i) For $\ep<0$, $p_\ep(\w,r)<0$ for all $\w\in\W$ if $r>s$ and $p_\ep(\w,r)>0$ for all $\w\in\W$ if $r<s$.
So, Remark \ref{rm:comparison}.1 shows that $\mA_\ep=\{s\}$.
\smallskip\par
(ii) Since $\{s\}$ is an attractive hyperbolic $\tau_0$-copy of the base,
Theorem \ref{th:cpositivo}(i) proves (ii).
\smallskip\par
(iii) For $\ep\in(0,\ep_*)$, $l_\ep<0$ determines an attractive hyperbolic copy of the base
(see Theorem \ref{th:cpositivo}(ii)), and $s$ another one. Hence, there exists a repulsive
hyperbolic copy of the base, $\{\mm_\ep\}$, with $\ml_\ep<\mm_\ep<s$: see Theorem \ref{th:nummin}.
For $\ep>\ep_*$, $\{\ml_\ep\}$ provides an attractive hyperbolic copy of the base and $\{s\}$ a repulsive one,
with $\ml_\ep<s$. Hence, Theorem \ref{th:nummin} ensures that $\{\muk_\ep\}$ is also an attractive
hyperbolic copy of the base above $\{s\}$.
\smallskip\par
(iv) Theorem \ref{th:cpositivo}(ii) shows that $\{\ml_{\ep_*}\}$ is an attractive hyperbolic
$\tau_{\ep_*}$-copy of the base. Since $\{s\}$ is a nonhyperbolic one, there are no more: see Theorem \ref{th:nummin}.
\smallskip\par
(v) Theorem \ref{th:cbpositivo}(ii) proves the assertions concerning $\ml_\ep$. Note now that the sets
$\mI$ of Theorem \ref{th:cpositivo}(iii) and $\mJ$ of Theorem \ref{th:cbpositivo}(iii) are $(0,\ep_*)$ and
$(\ep_*,\infty)$. Since $\muk_\ep>s$ for $\ep>\ep^*$, we get
$(\muk_{\ep_1})'(\w)\ge p_{\ep_2}(\w,\muk_{\ep_1}(\w))$ for all $\w\in\W$ if
$\ep_*<\ep_1<\ep_2$. According to Remark \ref{nota:Z}, this ensures that $\muk_{\ep_1}<\muk_{\ep_2}$. This
fact combined with
Theorem \ref{th:cbpositivo}(iii) proves the assertions concerning $\muk_\ep$ on $(\ep_*,\infty)$.
By reasoning as in the proof of Theorem \ref{th:cpositivo}(iii), we check that $(0,\ep_*)\to C(\W,\R),\,\ep\mapsto\mm_\ep$ is
strictly increasing, which combined with Theorem \ref{th:cpositivo}(iii) proves the assertions concerning $\muk_\ep$
on $[0,\ep_*)$.
The monotonicity properties ensure the existence of the limits $\lim_{\ep\to(\ep_*)^+}\muk_\ep\ge s$ and
$\lim_{\ep\to(\ep_*)^-}\mm_\ep\le s$, and that they provide semicontinuous $\tau_{\ep_*}$-equilibria, and
Proposition \ref{prop:defM} ensures that they coincide with $s$ at their continuity points.
\end{proof}
The bifurcation diagram described in Theorem \ref{th:caso3} is depicted in Figure \ref{fig:caso4}.

\begin{nota}\label{rm:extension}
For further purposes, we point out that the complete analysis can be repeated for
if we change \eqref{eq:EDOpol} by $x'=d(\w)\,p_\ep(\w,x)$ with
$d\colon\W\to(0,\infty)$ continuous and $p_\ep$ given by \eqref{def:p}.
In fact, $d$ does not has influence on the global shape of the bifurcation diagram
in the analysed cases: they still depend on the relation between $c$ and $-b/a$.
\end{nota}
We complete the description of these bifurcation diagrams by pointing out that
the bifurcation values $\ep_*$ and $\ep^*$ of Figure \ref{fig:caso1},
$\underline\ep$, $\bar\ep$ and $0$ (resp.~$0$) of the left (resp.~right) panel of
Figure \ref{fig:casos23}, and $0$ of Figure \ref{fig:caso4}, are points of discontinuity of
the map $\ep\to\mA_\ep$, in the sense explained in \cite[Chapter 3]{calr}. (In fact,
the map is upper-semicontinuous but not lower-semicontinuous at those points.)
\subsection{The results for an $\ep$-parametric family of ODEs}\label{subsec:process}
The results obtained so far in Section \ref{sec:results} provide a wealth of information about
the evolution of the dynamics induced by the $\ep$-parametric family of ODEs
\begin{equation}\label{eq:unaEDO}
 x'=\bar p_\ep(t,x)\,,
\end{equation}
where $\ep$ varies in $\R$ and
\[
 \bar p_\ep(t,x):=-x^3+\bar c(t)\,x^2+\ep\,\big(\bar b(t)\,x+\bar a(t)\big)
\]
for bounded and uniformly continuous maps $\bar c,\,\bar b,\,\bar a\colon\R\to\R$
such that the corresponding hull is minimal and uniquely ergodic.
\par
More precisely,
let $\W$ be the joint hull of $\bar\w:=(\bar c,\bar b,\bar a)$, defined as the closure in the compact
open topology of $C(\R,\R^3)$ of the time-shifts $\w_t=(\bar c_t,\bar b_t,\bar a_t)$ (see
Section \ref{subsec:hull}), and let us define $c(\w)=\w_1(0)$, $b(\w)=\w_2(0)$ and
$a(\w)=\w_3(0)$ for $\w=(\w_1,\w_2,\w_3)\in\W$. Note that $c(\bar\w_t)=\bar c(t)$,
$b(\bar\w_t)=\bar b(t)$ and $a(\bar\w_t)=\bar a(t)$.
That is, \eqref{eq:unaEDO} is one of the elements of the family
\begin{equation}\label{eq:unaEDOhull}
 x'=-x^3+c(\w_t)\,x^2+\ep\,\big(b(\w_t)\,x+a(\w_t)\big)\,,\qquad\w\in\W\,.
\end{equation}
There are well-known conditions ensuring the minimality and unique ergodicity of the
time-shift flow on $\W$, what we assume from now on.
For instance, this is the case if $\bar c,\,\bar b,\,\bar a\colon\R\to\R$
are almost periodic functions.

Let us represent by \eqref{eq:unaEDO}$_\ep$ the equation corresponding
to the value $\ep$.
The analysis of the global dynamics induced by \eqref{eq:unaEDO}$_\ep$ requires the
analysis of the set $\mB_\ep$ of bounded solutions of \eqref{eq:unaEDO}$_\ep$, which
is closely related to the attractor $\mA_\ep$ of \eqref{eq:unaEDOhull}$_\ep$: if $\ml_\ep$ and $\muk_\ep$
are the lower and upper bounded equilibria, and
if $l_\ep(t)=\ml_\ep(\bar\wt)$ and $u_\ep(t)=\muk_\ep(\bar\wt)$
(for $\bar\w=(\bar c,\bar b,\bar a)$), then $\mB_\ep:=\{(t,x)\,|\;l_\ep(t)\le x\le u_\ep(t)\}$.
This is a consequence of Theorem \ref{th:atractor}.
The analysis also relies on the number and type of hyperbolic solutions. Recall that
a bounded solution $b(t)$ of \eqref{eq:unaEDO}$_\ep$ is {\em hyperbolic attractive\/}
(resp.~{\em hyperbolic repulsive}) if there exist $k\ge 1$ and $\gamma>0$ such that
$\exp\Big(\int_s^t (\bar p_\ep)_x(r,b(r))\,dr\Big)\le ke^{-\gamma(t-s)}$ whenever $t\ge s$
(resp.~$\exp\Big(\int_s^t (\bar p_\ep)_x(r,b(r))\,dr\Big)\leq ke^{\gamma (t-s)}$ whenever $t\le s$).
According to, e.g., \cite[Theorems 5.3 and 5.6]{duno4},
\eqref{eq:unaEDO}$_\ep$ has at most three hyperbolic solutions,
in which case $l_\ep$ and $u_\ep$ are hyperbolic attractive and the ``middle one" is repulsive.
In this case, the dynamics of \eqref{eq:unaEDO}$_\ep$ is completely determined by its
hyperbolic solutions: see, e.g., \cite[Theorem 5.6]{duno4}.

It is easy to establish conditions on $\bar c$, $\bar b$ and $\bar a$ which are equivalent to
the hypotheses on $c$, $b$ and $a$ required in Theorems \ref{th:cpositivo}, \ref{th:cbpositivo},
\ref{th:caso1}, \ref{th:caso2} and \ref{th:caso3}:
\begin{itemize}[leftmargin=10pt]
\item[-] the conditions $c>0$, $a<0$, $b\ge 0$ and $b>0$ hold if and only if
$\inf_{t\in\R}\bar c(t)>0$, $\sup_{t\in\R}\bar a(t)<0$, $\inf_{t\in\R}\bar b(t)\ge 0$
and $\inf_{t\in\R}\bar b(t)>0$, respectively;
\item[-] the conditions $b(\w)\,c_++a(\w)<0$, $b(\w)\,c_-+a(\w)>0$, $c_+<3\,c_-$ and
$c_+<3\,s_-$ for $s_-=\inf_{\w\in\W}(-a(\w)/b(\w))$ are equivalent to
$\sup_{t\in\R}\bar c(t)<\inf_{t\in\R}(-\bar a(t)/\bar b(t))$,
$\inf_{t\in\R}\bar c(t)>\sup_{t\in\R}(-\bar a(t)/\bar b(t))$,
$\sup_{t\in\R}\bar c(t)<3\inf_{t\in\R}\bar c(t)$ and
$\sup_{t\in\R}\bar c(t)<3\inf_{t\in\R}(-\bar a(t)/\bar b(t))$, respectively;
\item[-] and the equality $a=-s\,b$ for a constant $s\in\R$ holds if and only if $\bar a=-s\,\bar b$.
\end{itemize}

Let us give an example of how to apply the previous results to the analysis of
the parametric variation of \eqref{eq:unaEDO}. (Another one, more precise, will be given
in Section \ref{sec:model}.)
\begin{prop}\label{prop:ejemploODE}
Let us assume that $\inf_{t\in\R}\bar c(t)>0$, $\inf_{t\in\R}\bar b(t)>0$,
$\bar a(t)=-s\,\bar b(t)$ for a constant $s\in(0,\infty)$ and all $t\in\R$, and $\inf_{t\in\R}\bar c(t)>s$.
Then,
\begin{itemize}[leftmargin=20pt]
\item[\rm(i)] $l_0$ is an attractive hyperbolic solution of \eqref{eq:unaEDO}$_0$
    and $u_0$ is a nonhyperbolic solution, with $\inf_{t\in\R}\bar c(t)<u_0(t)<\sup_{t\in\R}\bar c(t)$. In
    addition, $u_0(t)-c(t)$ changes sign at the points of a strictly increasing two-sided sequence $(s_n)_{n\in\Z}$.
\item[\rm(ii)] For all $\ep>0$, \eqref{eq:unaEDO}$_\ep$ has three hyperbolic solutions $l_\ep<m_\ep<u_\ep$, with
    $l_\ep<0<m_\ep<s<u_\ep$, with $l_\ep$ and $u_\ep$ attractive and $m_\ep$ repulsive. In addition,
    $\lim_{\ep\to\infty}u_\ep(t)=\infty$, $\lim_{\ep\to\infty}l_\ep(t)=-\infty$, $\lim_{\ep\to\infty}m_\ep(t)=s$
    and $\lim_{\ep\to 0^+}(u_\ep(t)-u_0(t))=\lim_{\ep\to 0^+}l_\ep(t)=\lim_{\ep\to 0^+}m_\ep(t)=0$,
    all of them uniformly on $\R$, and the maps $(0,\infty)\to C(\R,\R),$ $\ep\mapsto-l_\ep,m_\ep,u_\ep$ are
    strictly increasing.
\item[\rm(iii)] For all $\ep<0$, $\mA_\ep\subset\W\times(s,\sup_{t\in\R}\bar c(t))$;
    the map $(\infty,0]\to\R,\,\ep\mapsto u_\ep$ is strictly increasing;
    and there exists $\ep_0\le 0$ such that, for $\ep\in(-\infty,\ep_0)$,
    there is a unique bounded solution (given by $l_\ep=u_\ep$),
    which is hyperbolic attractive. If, in addition, $\sup_{t\in\R}\bar c(t)<3\,s$, then $\ep_0=0$,
\item[\rm(iv)] There are three possibilities for $\ep<0$:
    \begin{itemize}[leftmargin=10pt]
    \item[-] the value $\ep_0$ of {\rm(iii)} is $0$. (This is the situation if $\sup_{t\in\R}\bar c(t)<3\,s$.)
    \item[-] There exist $\underline\ep<\bar\ep<0$ such that: \eqref{eq:unaEDO}$_\ep$ has three
    hyperbolic solutions for any $\ep\in(\underline\ep,\bar\ep)$; $l_\ep=u_\ep$ and it is hyperbolic attractive
    for $\ep\in(0,\infty)-[\,\underline\ep,\,\,\bar\ep\,]$; $l_{\underline\ep}$ is the unique hyperbolic solution
    of \eqref{eq:unaEDO}$_{\underline\ep}$, it is attractive, and it is uniformly separated from $u_{\underline\ep}$;
    $u_{\bar\ep}$ is the unique hyperbolic solution
    of \eqref{eq:unaEDO}$_{\bar\ep}$, it is attractive, and it is uniformly separated from $l_{\bar\ep}$;
    the maps $(-\infty,\bar\ep)\to C(\R,\R),\,\ep\to l_\ep$ and $(\underline\ep,0)\to C(\R,\R),\,\ep\to u_\ep$
    are strictly increasing; and the map $(\underline\ep,\bar\ep)\to C(\R,\R),\,\ep\to m_\ep$
    is strictly decreasing.
\item[-] There is a unique point $\ep_1<0$ such that \eqref{eq:unaEDO}$_{\ep_1}$ has bounded solutions but not
    hyperbolic ones. In this case, $\inf_{t\in\R}(u_{\ep_1}(t)-l_{\ep_1}(t))=0$.
\end{itemize}
\end{itemize}
\end{prop}
The proof relies on applying Theorem \ref{th:caso2} to the families \eqref{eq:unaEDOhull}$_\ep$ constructed
from \eqref{eq:unaEDOhull}$_\ep$. It also uses that: a hyperbolic copy of the base $\{\mb\}$ for
\eqref{eq:unaEDOhull}$_\ep$ determines a hyperbolic solution of each equation \eqref{eq:unaEDOhull}$_\ep^\w$,
given by $t\mapsto \mb(\wt)$ (see, e.g.,~\cite[Proposition 2.7]{duno4}); and that, since $\W$ is minimal,
if $\mM$ is a nonhyperbolic $\tau_\ep$-minimal set and $(\w,x)\in\mM$, then the solution $v_\ep(t,\w,x)$
of \eqref{eq:unaEDOhull}$_\ep^\w$ is not hyperbolic (see, e.g., \cite[Proposition 1.54]{duen2}).
We leave the (easy) details to the reader.

We complete this section pointing out that, in the conditions of the previous result,
the possibly complicated dynamics arising at $\underline\ep$ reads as: it is possible
that $u_{\underline\ep}$ is the pointwise limit of $u_\ep$ and of $m_\ep$ as $\ep\to(\underline\ep)^+$, but not sure.
This happens if the point $\bar\w=(\bar c,\bar b,\bar a)$ belongs to a residual subset of $\W$, impossible to determine a priori,
which in addition can have measure ($m$) $0$ or $1$. In many situations (as when the initial coefficients are constants or
periodic maps) this residual set is the whole $\W$. If not, we just know that $u_{\underline\ep}\ge\bar u_{\underline\ep}\ge
m_{\underline\ep}$, with $\bar u_{\underline\ep}(t):=\lim_{\ep\to(\underline\ep)^+}u_\ep(t)$ and $m_{\underline\ep}(t):=\lim_{\ep\to(\underline\ep)^+}m_\ep(t)$.
The situation is analogous at $\bar\ep$ with $l_{\bar\ep}$ and limits as $\ep\to(\bar\ep)^-$; but not necessarily the same,
since the residual set can be different.
\section{Numerical simulations in a population dynamics model}\label{sec:model}
In this section, we study a single species population model that undergoes quasiperiodic
fluctuations (see for example \cite{RM} and references therein where experimental evidence
of quasi-periodic behavior in population dynamics can be found). We take into account
the interplay between the Allee effect (see for example \cite{cobg}) and migration phenomena,
both affected by seasonality. The model is
\begin{equation}\label{eq:model}
x'=r(t)\,x^2\left(1-x/k(t)\right)+\ep\,b(t)\,(x-s)\,,
\end{equation}
with $\ep\ge 0$.
The value $\ep=0$ of the parameter corresponds to the absence of migration. The maps $r(t)$
and $k(t)$ are positively bounded from below: $r(t)$
represents the intrinsic growth rate, i.e., the growth rate in case of unlimited resources, and
the function $k(t)$ is closely related to the carrying capacity although not exactly equal
(as in the autonomous case): it measures the threshold below which the per capita population
growth rate $x'/x$ decreases.
Changing the $x$ factor of logistic models to $x^2$ is a common way to include
the weak Allee effect: the per capita population growth rate is reduced at lower
density (as the $x=0$ solution loses hyperbolicity).

The additive term $\ep\,b(t)\,(x-s)$ is related to migration. The map $\ep\,b(t)$, where $b$
is positively bounded from below, represents the (seasonality-dependent) intensity of migration,
while $s$ is a positive constant representing the threshold of population attractiveness:
there is immigration (population increases) if the population is sufficiently high,
and emigration if it is below $s$. The idea fits well with that of the Allee effect:
a sufficient number of individuals increases the chances of survival of the population.
In fact, as we will explain later, for small and large values of $\ep>0$ (or even for all $\ep>0$)
there appear two strictly positive hyperbolic solutions of
\eqref{eq:model}$_\ep$: a critical population $m_\ep(t)$ (repulsive) that provides a threshold below
which the population will die out, and a stable (attractive)
healthy population $u_\ep(t)$ above this threshold. This is the
usual situation under strong Allee effect. Note that $x=0$ is not
a solution of \eqref{eq:model}$_\ep$ if $\ep>0$, so that extinction in finite
time is possible: at the moment in which $x$ reaches $0$,
the population disappears and the model becomes meaningless.
(Also observe that a negative value of $\ep$ (for $b>0$)
changes the role of $s$: there would be immigration for a lower number of individuals,
and this makes sense too. Although we will focus on $\ep\ge 0$, the previously obtained
theory also provides conclusions for $\ep<0$.)

A variation of $\ep$ means a variation in the migratory intensity which may give rise to
different population dynamics. The coefficient maps $r$, $k$ and $b$ are chosen to
get a quasiperiodic map $(r,k,b)\colon\R\to\R^3$, so the flow on its hull is
minimal and uniquely ergodic. Note that \eqref{eq:model} is
$x'=(r(t)/k(t))\,(-x^3+k(t)\,x^2+\ep\,k(t)\,b(t)\,(x-s)/r(t))$. Since
$\inf_{t\in\R}r(t)/k(t)>0$, the hull extension provides one of the families considered
in Remark \ref{rm:extension}, and hence all the results of Section \ref{sec:results}
can be applied to describe the bifurcation diagram
of the $\ep$-parametric family of skew-product flows. In the line of Section \ref{subsec:process},
in what follows we just focus on the family of processes instead of flows.
It is easy to check that the conditions $\sup_{t\in\R}k(t)<s$, $\inf_{t\in\R}k(t)>s$
and $k(t)\equiv s$ lead us respectively to the situations of Theorems \ref{th:caso1},
\ref{th:caso2} and \ref{th:caso3}, under their most restrictive hypotheses.
So, in all these cases, we have already proved the previously mentioned existence of two strictly
positive hyperbolic solutions $m_\ep<u_\ep$ for small or large $\ep>0$,
where $m_\ep$ is repulsive and $u_\ep$ is attractive, and where $\lim_{\ep\to 0^+}m_\ep(t)=0$
uniformly on $\R$. The solution $m_\ep$ acts as a threshold for survival:
if, for an $\ep>0$, $x_\ep(0)<m_\ep(0)$ or, equivalently, if $x_\ep(t_0)<m_\ep(t_0)$ for any $t_0>0$, then
the population becomes extinct in finite time; whereas, if $x_\ep(0)>m_\ep(0)$,
the population eventually ``reaches" (i.e., ``approaches until being undistinguishable
from") the healthy steady population $u_\ep(t)$.
In addition, $\lim_{\ep\to\infty}u_\ep(t)=\infty$ uniformly on $\R$, and so the resources that the
stable population requires exceed the capacity of the environment if $\ep$ is sufficiently large:
the increase of $u_\ep$ is due to a massive influx of individuals, difficult to imagine
for any reasonable population. But the model makes perfect sense for a not too large
migratory intensity. In addition, for $\ep=0$, any initial number of individuals gives rise to
a population which eventually reaches the (hyperbolic attractive) stable population $u_0(t)$,
where $\inf_{t\in\R}k(t)<u_0(t)<\sup_{t\in\R}k(t)$ for all $t\in\R$, and where $u_0-k$
changes sign infinitely many times as $t$ increases. This proves the
aforementioned close relation between $k(t)$ and the steady population in the absence of migration.

First, let us assume $\inf_{t\in\R}k(t)>s$, which is the situation of
Proposition \ref{prop:ejemploODE}: roughly speaking, the threshold of population
attractiveness is below the carrying capacity in the absence of migration. So, everything works
properly: only an initially too low population dies out, since, for ($\ep$-relatively) small $x$,
emigration dominates over intrinsic growth. More precisely,
Proposition \ref{prop:ejemploODE} ensures the absence of strictly positive bifurcation
values of $\ep$: the strictly positive hyperbolic solutions $m_\ep$ and $u_\ep$ exist
for all $\ep>0$. However, their monotonicity properties with respect to $\ep$
(see Figure \ref{fig:casos23}), and their behavior as $\ep\to 0^+$ and as
$\ep\to\infty$, give rise to a critical value of $\ep_{x_0}$ if we fix an
initial number of individuals $x_\ep(0)=x_0<s$: this ensures emigration for small $t>0$,
and implies that the population $x_\ep(t)$ survives reaching the steady one
if and only if this emigration is not too intense. More precisely,
there exists $\ep_{x_0}>0$ such that $x_\ep(0)<m_\ep(0)$ if
$0<\ep<\ep_{x_0}$ and $x_\ep(0)>m_\ep(0)$ if $\ep>\ep_{x_0}$, and this means
survival if $\ep<\ep_{x_0}$ (and in the unstable situation
$\ep=\ep_{x_0}$, with $x_{\ep_{x_0}}=m_{\ep_{x_0}}$),
and extinction in finite time for $\ep>\ep_{x_0}$. If the initial population is $x_0\ge s$,
then it survives for all $\ep>0$.

Now, let us analyze the case $k\equiv s$, adapting the information of Theorem \ref{th:caso3}
(see also Figure \ref{fig:caso4}). In this case, the threshold of attractiveness
coincides with the (constant) carrying capacity in the absence of migration.
If $\ep_*=s^2/\bar b$ for $\bar b=\lim_{t\to\infty}(1/t)\int_0^tb(s)\,ds$, then $s$ is
a constant solution for all $\ep>0$, hyperbolic attractive if and only if $\ep<\ep_*$
and hyperbolic repulsive if and only if
$\ep>\ep_*$. As in the previous case, an initial population $x_\ep(0)=x_0<s$ only survives
while the migration intensity is low enough to yield $x_0<m_\ep(0)$. Now,
this situation ends for sure at a value $\ep_{x_0}\le\ep_*$ of the parameter, beyond which our
population is doomed to extinction. Is addition, again as in the previous case, the population survives
for any $\ep$ if $x_0\ge s$. In terms of hyperbolic solutions,
there are two strictly positive ones $m_\ep<u_\ep$ for all $\ep>0$, $\ep\ne\ep_*$,
with $u_\ep=s$ for $\ep<\ep_*$ and $m_\ep=s$ for $\ep>\ep^*$. They approach each other as
$\ep\to\ep_*$, and their limits at this point are non uniformly separated bounded solutions
(which may coincide, as in the autonomous and periodic cases). So, there is a
(nonautonomous) transcritical bifurcation point at $\ep_*$.

In the rest of this section, we analyze the case $\sup_{t\in\R}k(t)<s$,
which we illustrate below with the help of a specific example. Now, the threshold of
population attractiveness is above the (again roughly speaking) carrying capacity of the
environment without migration. For instance, a low population may loose attractiveness even
with sufficient resources, leading to emigration; and immigration may occur if there are nearby
patches occupied by the same species but with fewer resources (and this causes population stress).
The precise information about the variation as $\ep$ increases is provided by Theorem
\ref{th:caso1} (see also Figure \ref{fig:caso1}), that shows the existence of two bifurcation values
$\ep_*<\ep^*$ on $(0,\infty)$ for the $\ep$-parametric family of skew-products. They can also
be read as bifurcation values for our initial process \eqref{eq:model}.
Note that a population that reaches a number of individuals less than $s$ stays below
$s$ (since $s$ is an upper solution), and its survival is only possible (not sure)
if $\ep\le\ep_*$ is small (see Figure \ref{fig:caso1}). In other words,
a low population is subject to emigration, and even if the first term on
the right-hand side of \eqref{eq:model}$_\ep$ is positive,
it is not sufficient to prevent extinction if the migration intensity is relatively large.

Let us now analyze the situation when $\ep\in(0,\ep_*]$. If $0<x_\ep(0)=x_0<m_{\ep_*}(0):=\lim_{\ep\to(\ep_*)^-}m_\ep(0)$,
then the population gets extinct in finite time when the migratory intensity exceeds
a certain value $\ep_{x_0}\in(0,\ep_*)$ (with $x_0=m_{\ep_{x_0}}(0)$).
If $x_0\ge m_{\ep_*}(0)$, then the population survives for all $\ep\in[0,\ep_*]$. But the steady population
$u_\ep$ eventually reached, which is below $s$, decreases as
$\ep$ increases: even if the initial population is above $s$ and hence there is
immigration during a time, the lack of resources causes a decrease, so the population
eventually reaches $s$, there is emigration forever, and a more intense emigration
means a lower steady population.

For $\ep\in(\ep_*,\ep^*)$, no matter how large the initial population is,
extinction arrive in finite time: the intensity of emigration once the population is below $s$
causes a decrease which is not compensated by the intrinsic growth rate. Finally,
for $\ep\ge\ep^*$: if $0<x_0<m_\ep(0)$ (as for all $\ep\ge\ep^*$ if $x_0\le s$ ),
then the population gets extinct in finite time; and if $x_0>s$,
then population survives for $\ep\ge\ep^{x_0}$, where this second critical
value $\ep^{x_0}\ge\ep^*$ satisfies $m_{\ep^{x_0}}(0)\le x_0$ (which exists since $m_\ep(0)$ decreases to
$s$ as $\ep\to(\ep^*)^+$). So, a high enough initial population can compensate the
stress caused by immigration. But now,
the only factor that allows the population to survive is immigration,
and hence the model looses credibility for a high value of $\ep$.

In terms of hyperbolic solutions, $\ep_*$ and $\ep^*$ are points at which two
hyperbolic solutions $m_\ep$ and $u_\ep$ approach each other (as $\ep\to(\ep_*)^-$
and as $\ep\to(\ep^*)^+$) , giving rise to the local absence of hyperbolic solutions:
they are replaced by limits of the monotonic families $\{u_\ep\}$ and $\{m_\ep\}$,
which may globally coincide or not, but which are never uniformly separated.

\begin{figure}
\includegraphics[width=\textwidth]{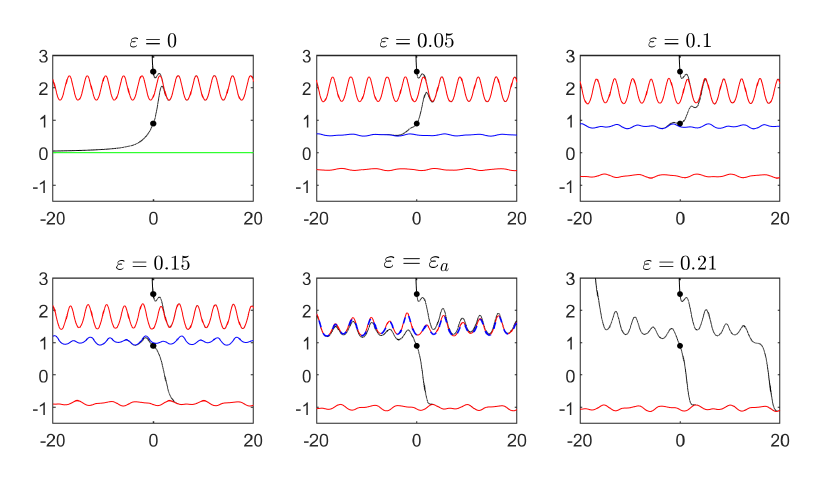}
\caption{Equation \eqref{eq:ex1}. Panels 1 to 5: Depiction of $u_{\ep}$,
$m_{\ep}$ and $l_{\ep}$ as $\ep$ increases towards $\ep_a$.
Panel 6: behavior of solutions for $\ep >\ep_a$.}\label{paper1}
\end{figure}
Let us illustrate these theoretical results with a particular example. The
quasiperiodic fluctuations in the population dynamics are represented by
taking $r\equiv 1$, and $k>0$ and $b>0$ periodic with incommensurate
oscillation frequencies. More precisely, we work with
\begin{equation}\label{eq:ex1}
 x'=-\frac{1}{2+0.5\sin(\sqrt 3\,t)}\:x^3+x^2+\ep\,(2.1+0.3\cos(t))\,(x-2.6)\,,
\end{equation}
so that $b(t):=2.1+0.3\cos(t)$ and $s:=2.6>2+0.5\sin(\sqrt 3\,t)=:k(t)$ for all $t\in\R$: we are in the
third of the cases described above. In what follows, we
detect numerically the bifurcation values $\ep_*$ and $\ep^*$ for \eqref{eq:ex1}
and approximate the hyperbolic copies of the base for different parameter values.

In order to approximate $\ep_*$, we reason as follows.
We choose $\ep_a^0$ and $\ep_b^0$ such that \eqref{eq:ex1}$_\ep$ has three hyperbolic solutions
for $\ep=\ep_a^0$ and just one for $\ep=\ep_b^0$. The existence of the hyperbolic solutions
has been detected numerically solving \eqref{eq:ex1} with a $4$-th order Runge-Kutta method
and constant discretization stepsize $h=2^{-10}$.  We obtain analogous numerical results also
for smaller and larger stepsizes. We solve initial value problems in $[-10^4, 10^4]$ for
forward integration and in $[10^4, -10^4]$ for backward integration.
We then apply a bisection procedure to the starting interval $[\ep_a^0, \ep_b^0]$ and
locate $\ep_*$ in $[\ep_a,\ep_b]=[0.201945926862769 \,\ 0.201945926863700]$. Note that
$(\ep_b-\ep_a)=O(10^{-12})$.
We reason in a similar way for $\ep^*$ and locate it inside the interval
$[9.129175817935083, \,\ 9.129175817935174]$.

\begin{figure}
\begin{minipage}[c]{0.48\linewidth}
\centering
\includegraphics[width=6cm]{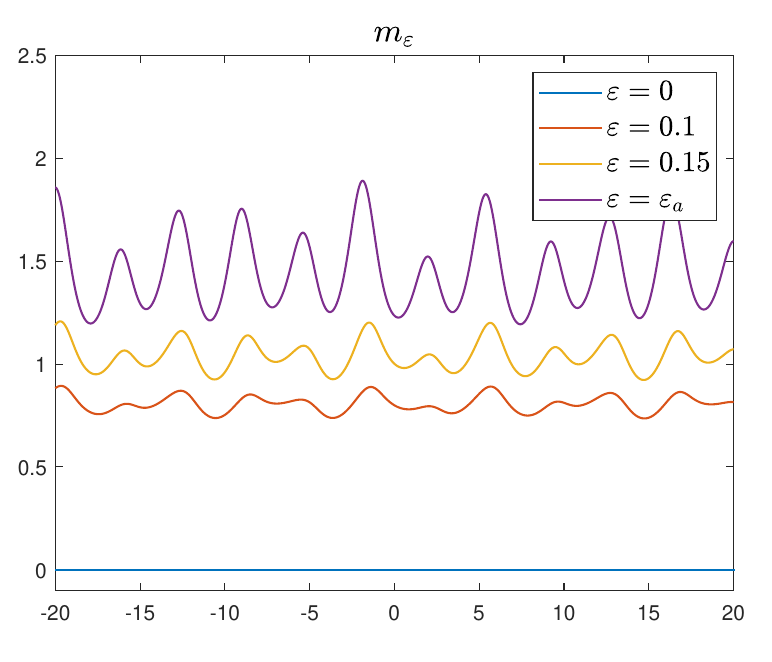}
\end{minipage}
\begin{minipage}[c]{0.48\linewidth}
\centering
\includegraphics[width=6cm]{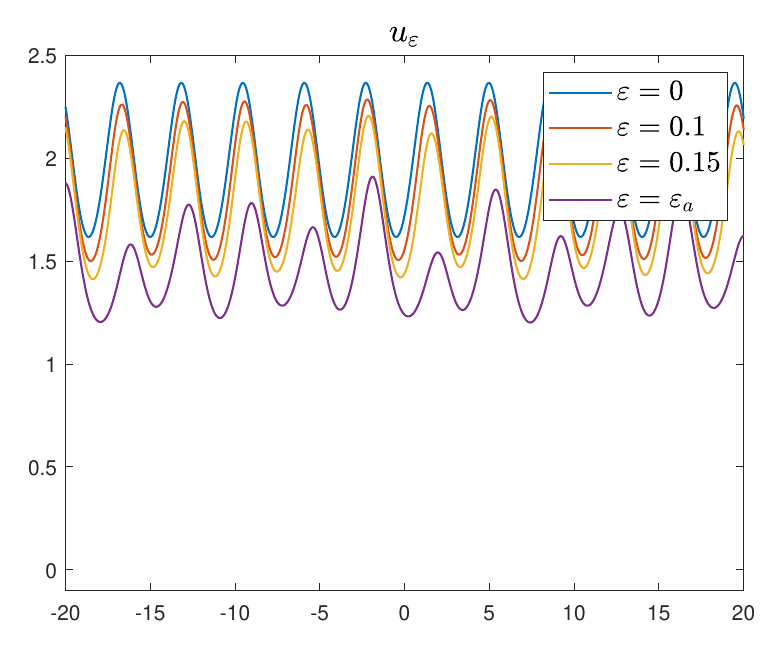}
\end{minipage}
\caption{Equation \eqref{eq:ex1}. Behavior of $m_{\ep}$ and $u_{\ep}$ as
$\ep$ increases towards $\ep_a$. The maps $m_{\ep_a}$ and $u_{\ep_a}$ are
indistinguishable in the representation window.}\label{paper3}
\end{figure}

In Figure \ref{paper1} we depict solutions of \eqref{eq:ex1} for different values of
$\ep$ in a neighborhood of $\ep_*$. In all six panels we plot in red the attractive hyperbolic
solutions (including the lower bounded solution $l_\ep$, which is negative for $\ep>0$),
in blue the repulsive one, and in black other solutions. We use same initial conditions
in all panels for the solutions in black, namely $x_\ep(0)=2.5$ and $x_\ep(0)=0.9$.
These initial data are marked on all panels.
The first panel of Figure \ref{paper1} corresponds to $\ep=0$: the model
does not contemplate migration. In addition to the already described behavior,
it is worth to observe that an initially low number of individuals results in
a low population over a long period (as a consequence of the weak Allee effect).
The solution $l_0=m_0=0$ is depicted in green, since it is non hyperbolic.
Panels $2,3,$ and $4$ correspond to three values of $\ep<\ep_*$, so that \eqref{eq:ex1}$_\ep$ has
three hyperbolic solutions: $u_{\ep}$ and $l_{\ep}$, depicted in red and
$m_{\ep}$, in blue. The solution $m_\ep$ acts as a threshold for survival.
We call attention to the variation of the solution depicted in black with initial
condition $x_\ep(0)=x_0=0.9$: in panels $2$ and $3$, it converges towards
$u_{\ep}$, while in panel $4$ it converges towards $l_{\ep}$ and hence
it gets extinct in finite time; so, the critical value $\ep_{x_0}$ before described
is inside $(0.1,0.15)$. For $\ep=\ep_a$ we can not distinguish between $u_{\ep}$ and $m_{\ep}$:
the two solutions seem to collide (for sure, they are not uniformly separated at $\ep_*$)
and hence they loose hyperbolicity. Populations above $u_{\ep_a}$ survive,
while populations below it get extinct. Finally for $\ep=0.21$, which is between
$\ep_*$ and $\ep^*$, the population is always doomed to extinction.

\begin{wrapfigure}{r}{0.5\textwidth}
\begin{center} \includegraphics[width=0.45\textwidth]{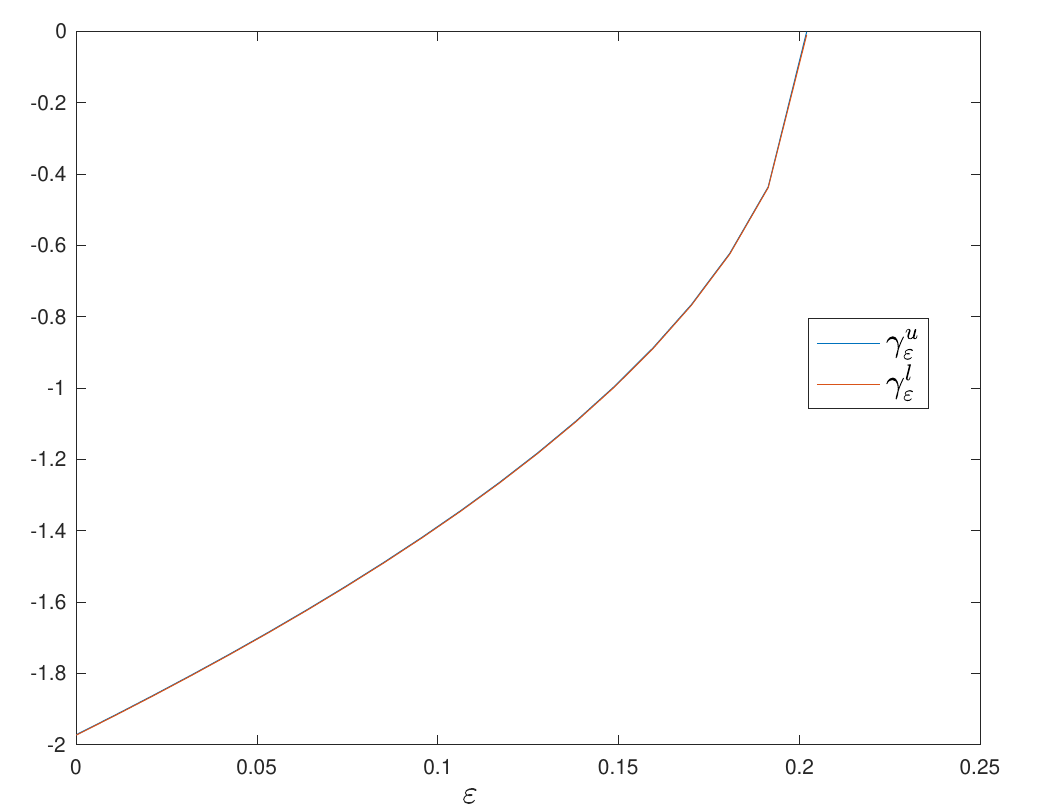}
\caption{Truncated Lyapunov exponents of $u_\ep$ as $\ep$ varies.}
\label{fig:Lyap}
\end{center}
\end{wrapfigure}
To complete the analysis, we explore the behavior of $u_{\ep}$ and $m_{\ep}$ as $\ep$ increase towards
$\ep_a$, to graphically show that they approximate each other as well as their loss of
hyperbolicity as $\ep\to(\ep_*)^+$. In the first and second panel of
Figure \ref{paper3} we respectively plot $u_{\ep}$ and $m_{\ep}$ for different
values of $\ep$. The plots show that $u_{\ep}$ decreases towards $u_{\ep_a}$ and
that $m_{\ep}$  increases towards $m_{\ep_a}$: see Theorem \ref{th:caso1}(i).
In order to verify the loss of hyperbolicity, we compute the Lyapunov exponent of
$u_{\ep}$ as $\ep$ increases towards $\ep_a$. The Lyapunov exponent is computed
truncating \eqref{def:2gamma} at a large enough time $T$. Equation \eqref{eq:ex1}
is quasiperiodic, hence the Lyapunov exponent of $u_{\ep_a}$ exists as a limit.
Nonetheless, we compute upper and lower approximations of $\gamma_{\ep_a}$ in order
to locate the Lyapunov exponent in a given interval. To this purpose, together with $T$,
we also use a finite time $\tau$, with $0 \ll \tau<T$. We then take as lower (resp. upper)
approximation the minimum (maximum) of all the truncated exponents for $t \in [\tau, T]$.
We denote this lower and upper approximations respectively as $\gamma^l_{\ep_a}$ and
$\gamma^u_{\ep_a}$. The values that we obtain are showed in Table \ref{table:les}.
From this table, the linear convergence of the Lyapunov exponent to $0$ is evident,
confirming the loss of hyperbolicity of $u_{\ep}$ at the bifurcation value.
Finally, in Figure \ref{fig:Lyap}, we plot the Lyapunov exponents $\gamma^l_{\ep}$ and
$\gamma^u_{\ep}$  of $u_{\ep}$ as $\ep$ varies in $[0, \ep_a]$. In the plot, the upper and lower
exponent approach zero as $\ep$ approaches $\ep_a$, witnessing a loss of hyperbolicity of
$u_{\ep}$ as $\ep \to \ep_*$.

\begin{table}[H]
\begin{tabular}{|c|c|c|c|}
\hline
T & $\tau$ & $\gamma^l_{\ep_a}$ & $\gamma^u_{\ep_a}$ \\
\hline
$10^3$ & $10^2$ & $-1.8 \times 10^{-2}$ & $\quad 3.1 \times 10^{-3}$ \\
\hline
$10^4$ & $10^3$ & $-9.4 \times 10^{-3}$ & $-1.2 \times 10^{-3}$ \\
\hline
$10^5$ & $10^4$ & $-1.3 \times 10^{-3}$ & $-1.5 \times 10^{-4}$ \\
\hline
\end{tabular}
\vspace{8pt}
\caption{Lower and upper approximations of the Lyapunov exponent of
$u_{\ep_a}$ for different values of the truncated time $\tau$ and $T$.}\label{table:les}
\end{table}


\end{document}